\documentclass[11pt]{article}
\usepackage{pstricks-add}
\usepackage{pgfplots}
\usepackage[latin2]{inputenc}
\usepackage{multicol}
\usepackage{t1enc}
\usepackage{babel}
\usepackage{enumerate}
\usepackage{graphicx}
\usepackage{float}
\usepackage{amstext}
\usepackage{amsthm}
\usepackage{amssymb}
\usepackage{amsmath}
\usepackage{color}
\usepackage[a4paper]{geometry}
\usepackage{pst-plot}
\usepackage{pst-eucl}
\pgfplotsset{compat=1.16}
\usetikzlibrary{shadows.blur}
\usepackage{tikz}
\usepackage{hyperref}
\usepackage{mathtools}
\usepackage{tkz-euclide}
\usepackage{fancyhdr}

\usepackage{hyperref}
\hypersetup{
	colorlinks=true,
	linkcolor=blue,
	filecolor=magenta,      
	urlcolor=cyan,
	pdftitle={Overleaf Example},
	pdfpagemode=FullScreen,
}

\hypersetup{
	colorlinks=true,
	linkcolor=black,
	filecolor=black,      
	urlcolor=black,
	pdftitle={Sharelatex Example},
	bookmarks=true,
}

\newcommand\p{\mathbf P}

\newcommand\R{\mathbb{R}}
\newcommand\N{\mathbb{N}}

\newcommand\dd{\mathrm{d}}

\newcommand{\ans}[1]{\noindent\textbf{Solution. }#1}

\newcommand{\sna}[1]{\noindent\textbf{Remark. }#1}

\newtheorem{theorem}{Theorem}
\theoremstyle{definition}
\newtheorem{fa}{Exercise}

\theoremstyle{definition}

\newtheorem*{fel*}{Feladat}

\newtheorem*{meg*}{Solution}
\newtheorem*{megv*}{Megoldás vázlat}
\newtheorem*{megj*}{Megjegyzés}
\newtheorem*{def*}{Definition}

\makeatletter
\tikzset{arrow coord style/.style={%
		densely dashed,
		\tkz@euc@linecolor,
}}
\tikzset{xcoord style/.style={%
		\tkz@euc@labelcolor,
		font=\normalsize,text height=1ex,
		inner sep = 0pt,
		outer sep = 0pt,
		fill=\tkz@fillcolor,
		below=6pt
}} 
\tikzset{ycoord style/.style={%
		\tkz@euc@labelcolor,
		font=\normalsize,text height=1ex, 
		inner sep = 0pt,
		outer sep = 0pt, 
		fill=\tkz@fillcolor,
		left=6pt
}}  
\makeatother

\begin{document}
	
\title{The two-sided approximation and counting method}
\author{\textsc{Attila M\'ader}\thanks{
			Bolyai Institute, University of Szeged, 
			Aradi v\'ertan\'uk tere 1, 6720 Szeged, Hungary; 
			e-mail: \texttt{madera@math.u-szeged.hu}}, \textsc{M\'at\'e Szalai}\thanks{
			Bolyai Institute, University of Szeged, 
			Aradi v\'ertan\'uk tere 1, 6720 Szeged, Hungary; 
			e-mail: \texttt{szalaim@math.u-szeged.hu}}}
\maketitle

\begin{abstract}
Below, we summarize the appearances and possible uses of the two-sided approach and the two-sided counting in the most diverse areas of (secondary) school mathematics.

\vspace{0.4cm}

\noindent \textit{Keywords}: Approximation, model description, two-sided counting.

\vspace{0.4cm}

\noindent \textit{MSC2020}: 97D20, 97N50.

\end{abstract}

\section{Introduction}

One of the topics in the oral part of advanced mathematics is regularly called \textit{Methods of Proof and Their Presentation in the Proof of Theorems}. It is a problem not only for the students, but sometimes even for the preparatory or examiners, exactly which methods can be considered as belonging to the requirements of the core material, and which methods should be included in the short time for the exam, through which items, the candidate to keep it short, understandable, but not trivial. Most aids mention the following: direct proof, indirect proof, induction, Pigeonhole principle.

The method of two-sided addition is nowhere, or very rarely, included. Nevertheless, this method, and in general the two-sided approach, can be used in many areas of mathematics to avoid and simplify formal proofs, as well as to ,,nice'', essential, illustrative proofs. The method was first used by Archimedes (287 BC - 212 BC). Changing points of view is also particularly difficult in mathematics.

\begin{figure*}[h!]
	\centering
	\includegraphics[width=0.3\linewidth]{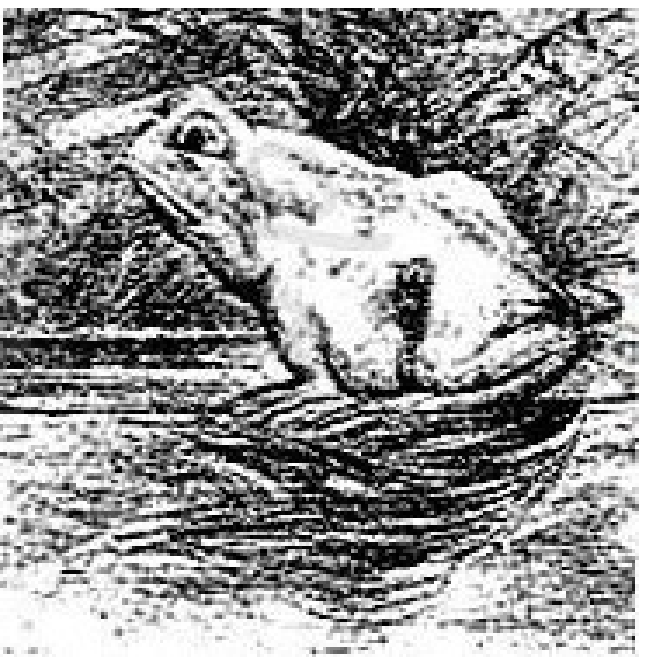}\hskip2cm\includegraphics[width=0.3\linewidth]{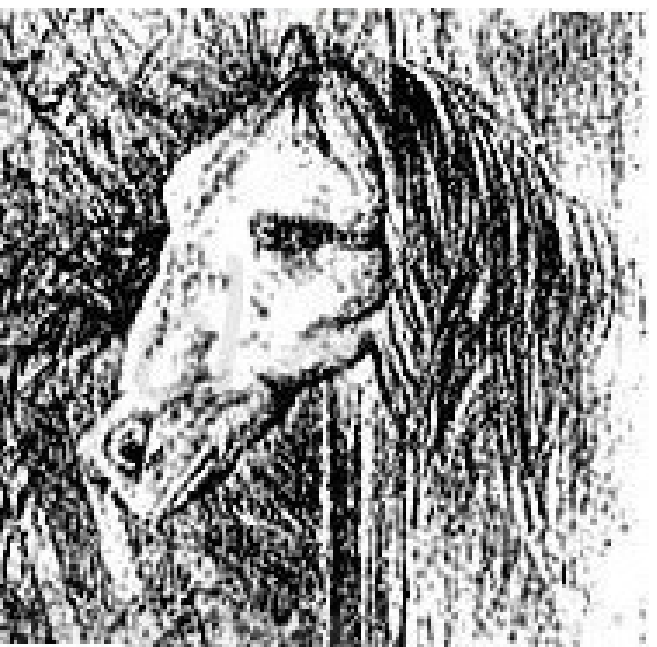}
	\caption{Frog or horse? - both a frog and a horse}
	\label{fig:beka}
\end{figure*}

However, the multiple approaches, to the same thing from several perspectives, directions, methods allow for behind-the-scenes, exploration of real causes, access to actual mathematical content (Figure \ref{fig:beka}).

\section{Algebra}

\subsection{Operational properties}

\begin{figure}[h!]
	\centering
	\includegraphics[width=0.3\linewidth]{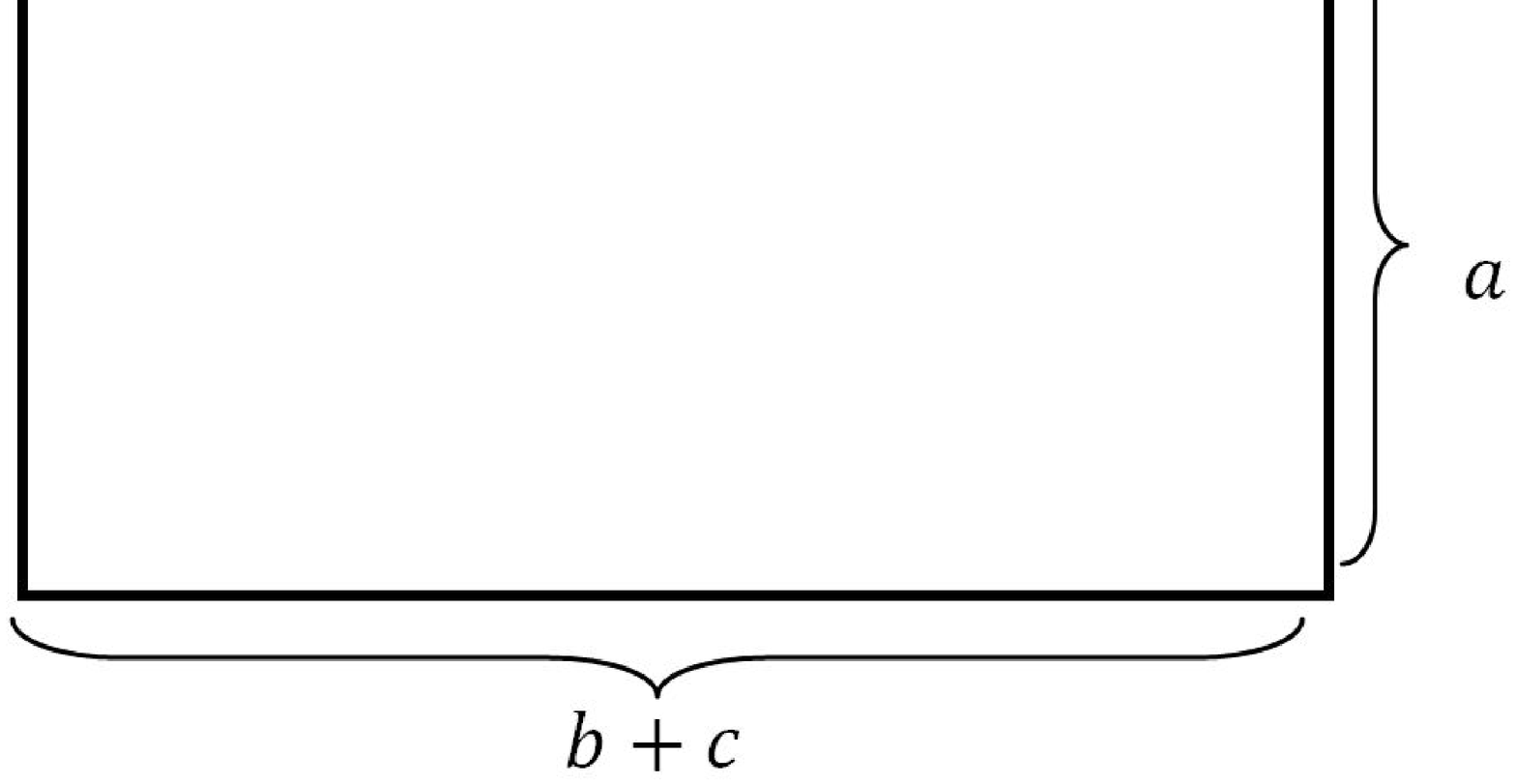}\includegraphics[width=0.3\linewidth]{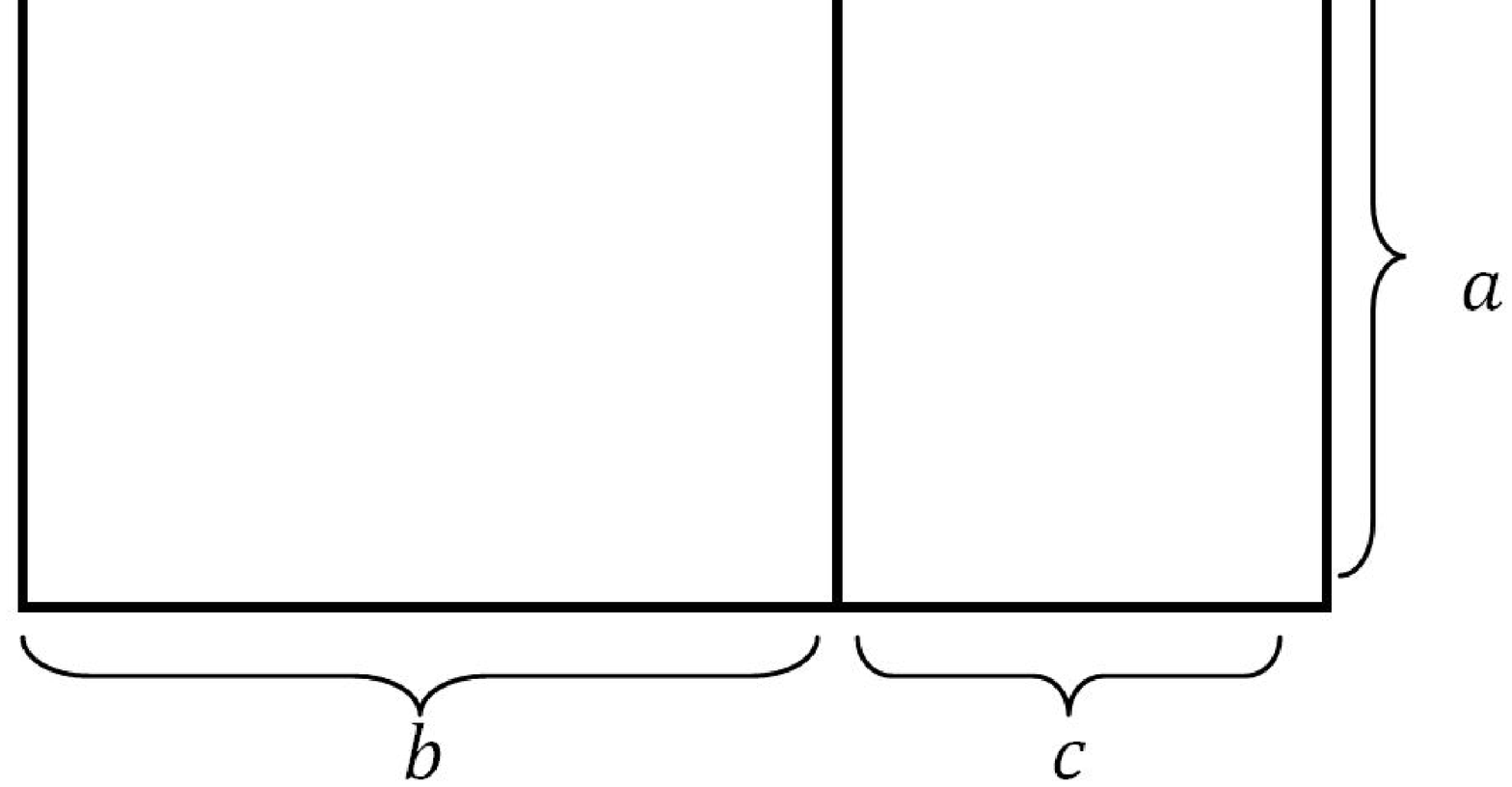}
	\vskip-4cm
	\caption{$(b+c)a=ba+ca$}
	\label{fig:k01}
\end{figure}

In elementary school, algebraic expressions first appear in their abstract form. At this age, students' abstraction skills are not yet at such a level that they can master $(b+c)a=ba+ca$ type identities (operational properties) in an abstract way. After the learning difficulties and often failures of identities that are often intended to be acquired only in abstract form, it is no wonder that many people have trouble dealing with similar expressions later on. However, by using a simple geometric model - the (so far) known method of calculating the area of a rectangle-, we can help the understanding with an diagram, for example as follows (Figure \ref{fig:k01}). 

The figure on the left shows a rectangle with sides $b+c$ and $a$, respectively, so its area is $(b+c)a.$ Cut our rectangle into two rectangles as shown in the figure on the right. Thus we obtained two rectangles with sides $b$ and $a$ and $c$ and $a$. Of course, the sum of the areas of the two parts is equal to the area of the original rectangle, i.e. $(b+c)a=ba+ca.$

\sna{In the product of the sum of the two terms by the sum of the two terms, the formalism of \textit{,,every term with every term''} can be made easy to understand and can be recalled later with the help of a visual representation similar to the above. Now, with the help of a small trick, we will cut a suitable rectangle with cuts parallel to its sides, and write down its area before and after cutting. Consider the following figure (Figure \ref{fig:k04}). 
	
\begin{figure}[h!]
	\centering
	
	\begin{minipage}{.5\textwidth}
		\centering
		\includegraphics[width=0.7\linewidth]{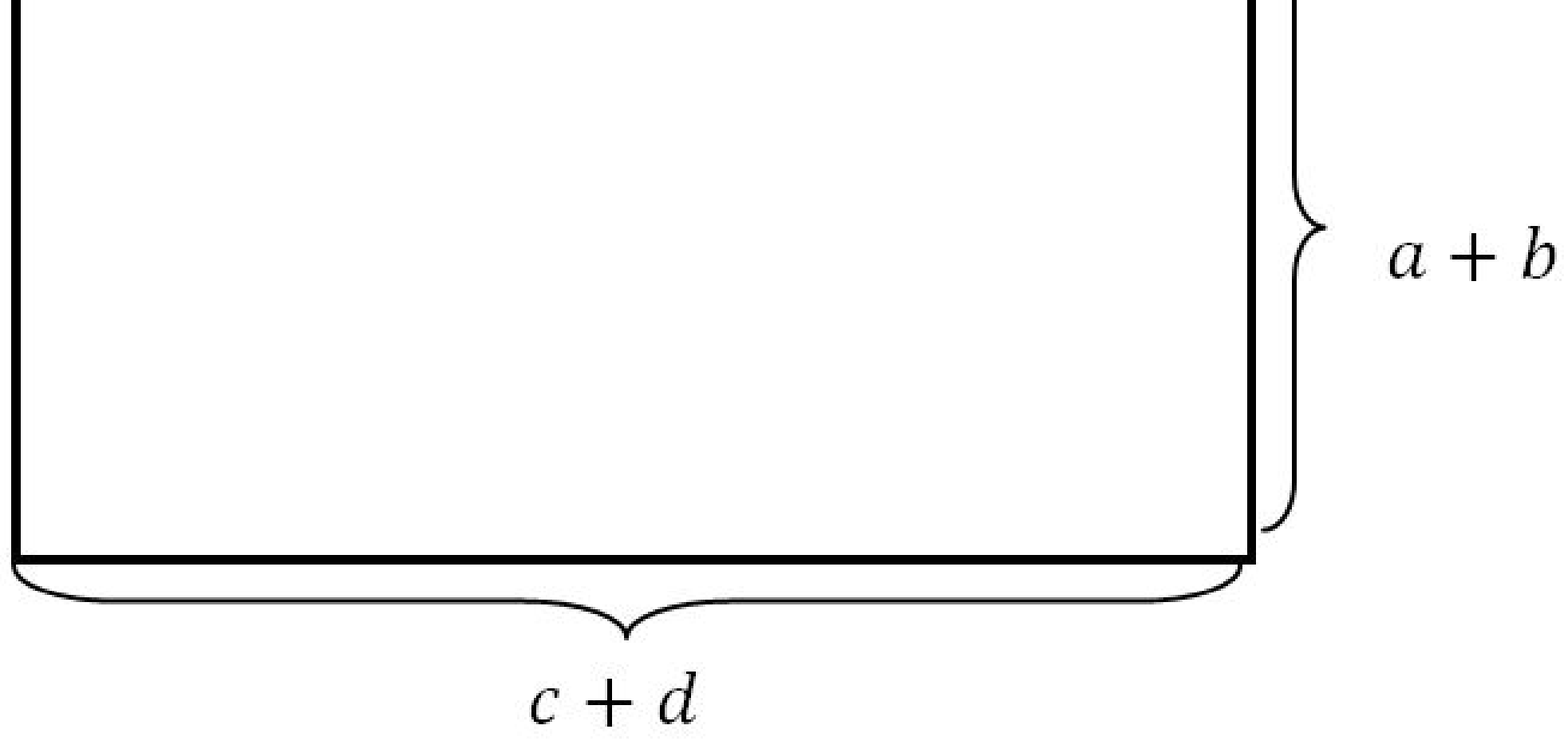}
		\label{fig:test1}
	\end{minipage}%
	\hspace{-3cm}
	\begin{minipage}{.5\textwidth}
		\vspace{-1.3cm}
		\centering
		\includegraphics[width=0.6\linewidth]{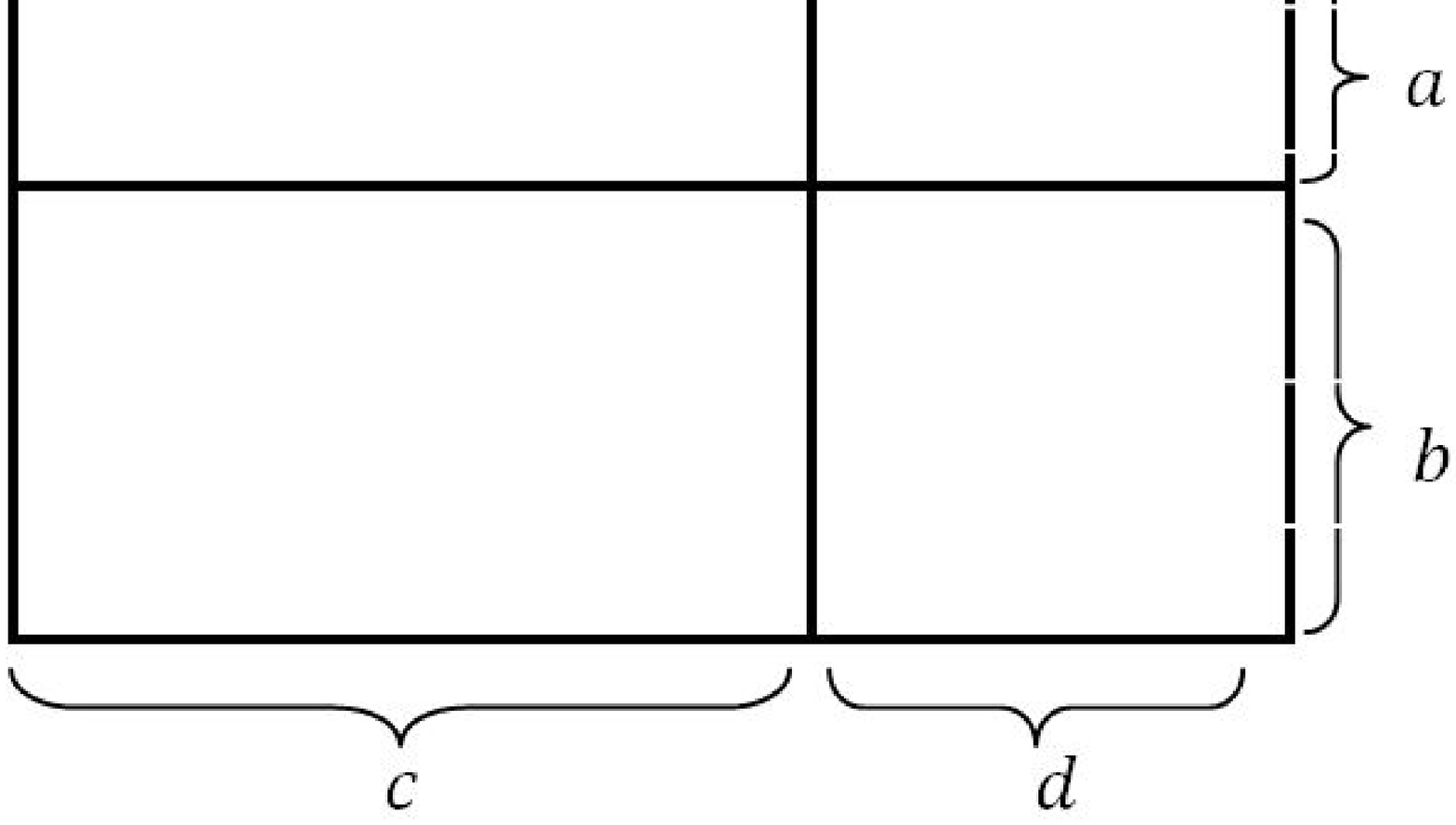}
		\label{fig:test2}
	\end{minipage}
	\vskip-4.2cm
	\caption{$(a+b)(c+d)=ac+ad+bc+bd$}
	\label{fig:k04}
\end{figure}}

\subsection{Notable identities}

A significant part of the ninth grade curriculum deals with notable multiplications, different methods of multiplication, and operations with algebraic fractions. Possessing this knowledge at the asset level is a prerequisite for progress, but learning it is difficult. One reason for this is that their teaching, including applications, is formal.

\begin{figure}[h!]
	\centering
	\includegraphics[width=0.3\linewidth]{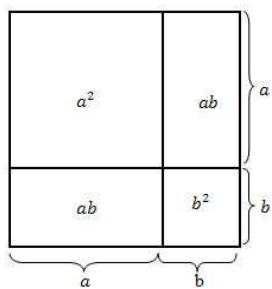}
	\vskip-2cm
	\caption{$(a+b)^2=a^2+2ab+b^2$}
	\label{fig:k07}
\end{figure}

By modifying the above figures that can be used in primary school (e.g. if the corresponding two rectangles in Figure \ref{fig:k04} are square) we can get help that can be used in high school as well (Figure \ref {fig:k07}). By writing the area of the square in two ways, you get: $(a+b)^2=a^2+2ab+b^2$.

\begin{fa}
	
	Let us prove and illustrate that 
	
	\begin{multicols}{2}
	\begin{enumerate}[(a)]
		\item $a(b-c)=ab-ac,$
		\item $(a-b)(c+d)=ac+ad-bc-bd,$
		\item $(a-b)(c-d)=ac-ad-bc+bd,$
		\item $(a-b)^2=a^2-2ab+b^2.$
		\item $(a+b)^3=a^3+3a^2b+3ab^2+b^3.$
	\end{enumerate}
	\end{multicols} 

\end{fa}

\subsection{Notable amounts - $\mathbf{1+2+3+\ldots+n}$}

\begin{figure}[h!]
\centering
\includegraphics[width=0.2\linewidth]{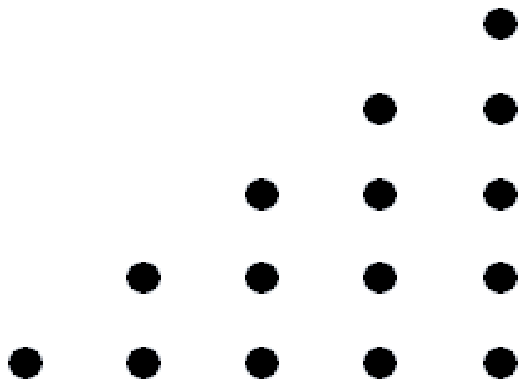}\includegraphics[width=0.2\linewidth]{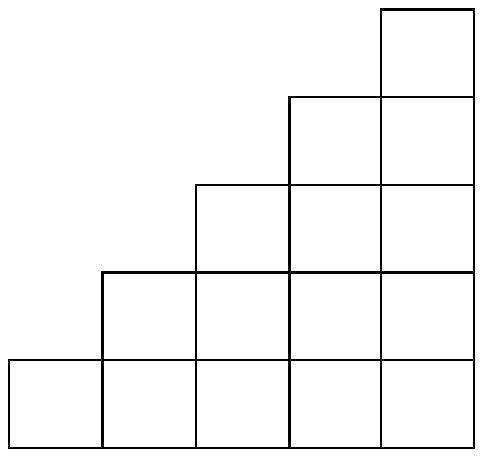}
\caption{$1+2+3+\ldots$}
\label{fig:k1}
\end{figure}

Many people know Gauss's beautiful proof of the closed form of the sum of the first $n$ positive integers based on oral tradition. The relation thus obtained:

\begin{equation*}\label{rt}
1+2+3+\ldots+n=\frac{n(n+1)}{2}
\end{equation*} 
studied, despite the evidence, the question of \textit{why}, precisely \textit{why exactly this}, may often arise in our students. The question can be answered, and thus the essence of the relationship, the true mathematical content behind the formula can be illuminated by carrying out a suitably chosen enumeration in two ways. Of course, all this can be helped with a good diagram. The sum of the first $n$ positive integers can be naturally illustrated as follows. Place pebbles (squares of unit area) as shown in the figure \ref{fig:k1}.

\begin{figure}[h!]
\centering
\includegraphics[width=0.2\linewidth]{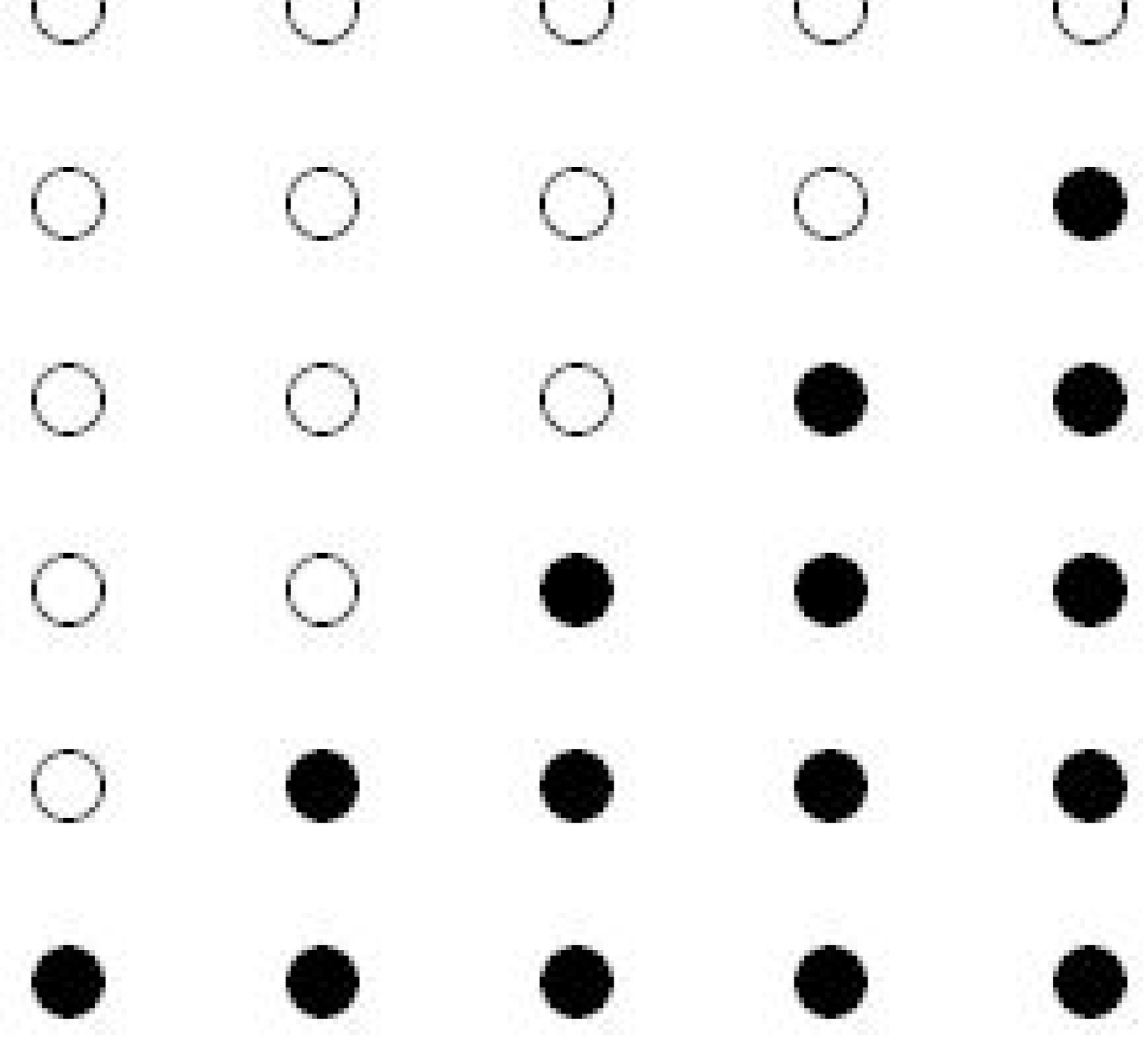}\includegraphics[width=0.2\linewidth]{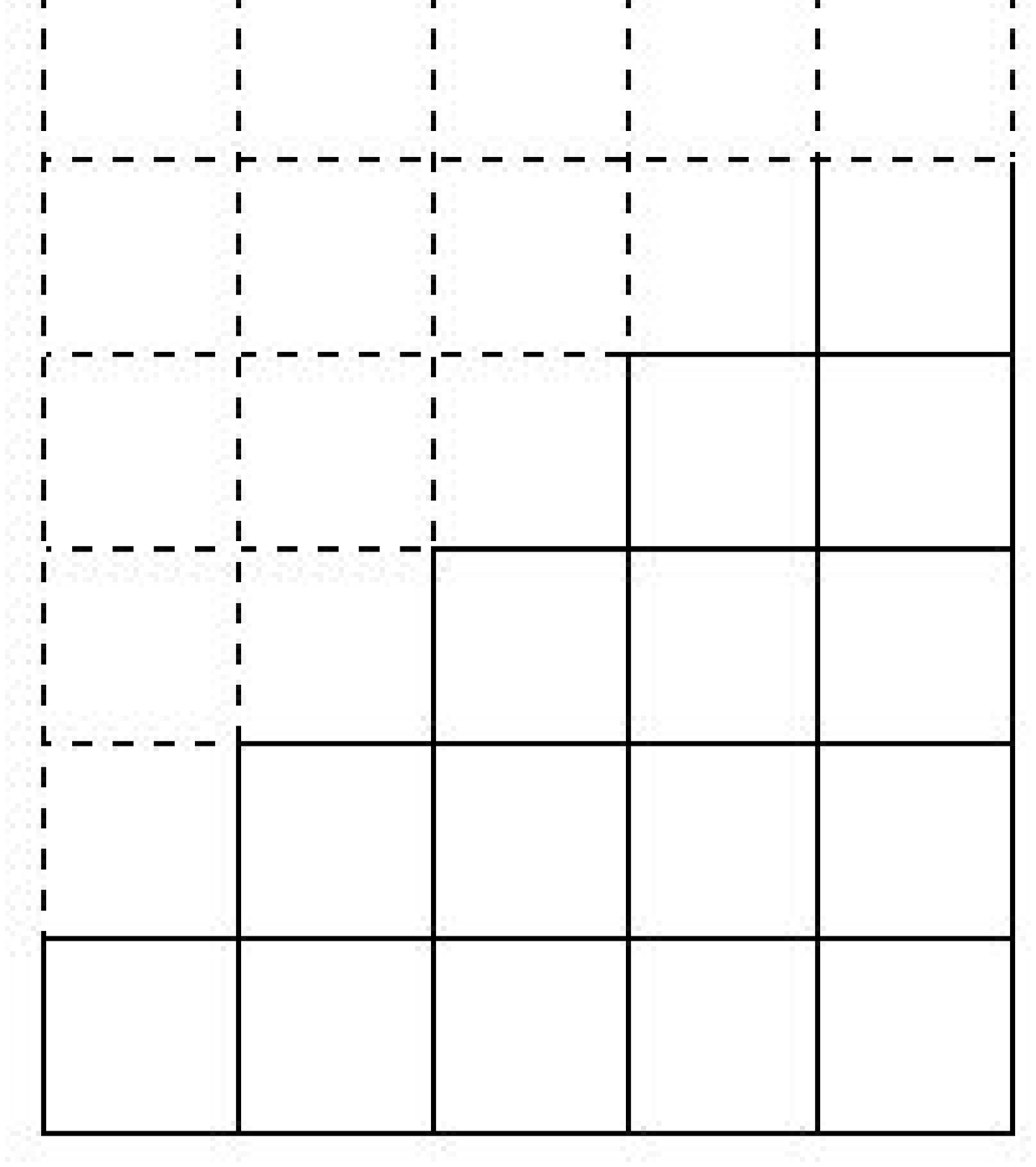}
\vskip-1.5cm
\caption{The $n\times (n+1)$ rectangle and its half}
\label{fig:k08}
\end{figure}

What we got is half of a $n\times (n+1)$ rectangle (Figure \ref{fig:k08}).

The number of pebbles in the above arrangement can be determined in two ways. Summing the pebbles (area of squares) per column, the number of pebbles (the area): $1+2+\ldots +n.$ By counting the pebbles (areas of squares) at the same time, using the half rectangle, we get $\frac{n(n+1)}{2}$ number of pebbles. Thus, by counting the number of pebbles in two ways, we get: $$1+2+\ldots+n=\frac{n(n+1)}{2}.$$

\begin{fa}
	\begin{enumerate}[(a)]
		\item The elements of the series $1, 1+2, 1+2+3, \ldots$ are also called triangular numbers. Why?
		\item We know that the difference between consecutive triangular numbers is respectively $2, 3, \ldots.$ But what can we say about the sum of two adjacent triangular numbers?
	\end{enumerate}
\end{fa}

\subsection{Notable sums - $\mathbf{1+3+5+\ldots+2n-1}$}

The fact that the sum of the first $n$ positive odd numbers is exactly $n^2$ is particularly surprising. Of course, the proof can be done by applying the versatile method, induction, but there is definitely a proof suitable for the ,,Big Book''\footnote{According to Pál Erdős, there must be a ,,Big Book'', with the most beautiful proof of all theorems.} for a nice connection. Maybe the next one is just like that. The search for a new method is necessary not only because of the avoidance of empty formalism, but also because the induction proof in this case does not answer \textit{why}, even though we can definitely expect this from a good proof.

Place $n^2$ number of pebbles in the shape of a $n\times n$ square (Figure \ref{fig:k13}). Divide our square, as shown in the figure. From determining the number of pebbles in two ways: $$1+3+\ldots+(2n-1)=n^2.$$

\begin{figure}[h!]
	\centering
	\definecolor{uuuuuu}{rgb}{0.26666666666666666,0.26666666666666666,0.26666666666666666}
	\begin{tikzpicture}[scale=0.65][line cap=round,line join=round,>=triangle 45,x=1.0cm,y=1.0cm]
		\clip(-0.5,-0.5) rectangle (11.,5.);
		\draw [line width=2.pt] (5.60345493593774,4.420107167630347)-- (10.388268164892615,4.3887655962616465);
		\draw [line width=2.pt] (5.572113364569041,3.4067296933757136)-- (9.40623226200668,3.385835312463247);
		\draw [line width=2.pt] (5.593007745481507,2.414246600033547)-- (8.392854787752045,2.3829050286648474);
		\draw [line width=2.pt] (5.60345493593774,1.4113163162351476)-- (7.389924503953644,1.3799747448664477);
		\draw [line width=2.pt] (5.593007745481507,0.4083860324367478)-- (6.407888601067708,0.3979388419805145);
		\draw [line width=2.pt] (6.397175294831571,0.39807619206046524)-- (6.386994220155242,-0.1975510140247854);
		\draw [line width=2.pt] (7.389924503953644,1.3799747448664477)-- (7.41081888486611,-0.23933977584971872);
		\draw [line width=2.pt] (8.392854787752045,2.3829050286648474)-- (8.382407597295812,-0.21844539493725204);
		\draw [line width=2.pt] (9.385395432825325,3.3859488646658704)-- (9.385337881094214,-0.22889258539348536);
		\draw [line width=2.pt] (10.388268164892615,4.3887655962616465)-- (10.409162545805081,-0.1975510140247854);
		\begin{scriptsize}
			\draw [fill=uuuuuu] (0.,0.) circle (4.5pt);
			\draw [fill=black] (1.,0.) circle (4.5pt);
			\draw [fill=black] (2.,0.) circle (4.5pt);
			\draw [fill=black] (3.,0.) circle (4.5pt);
			\draw [fill=black] (4.,0.) circle (4.5pt);
			\draw [fill=black] (0.,1.) circle (4.5pt);
			\draw [fill=black] (1.,1.) circle (4.5pt);
			\draw [fill=black] (2.,1.) circle (4.5pt);
			\draw [fill=black] (3.,1.) circle (4.5pt);
			\draw [fill=black] (4.,1.) circle (4.5pt);
			\draw [fill=black] (0.,2.) circle (4.5pt);
			\draw [fill=black] (1.,2.) circle (4.5pt);
			\draw [fill=black] (2.,2.) circle (4.5pt);
			\draw [fill=black] (3.,2.) circle (4.5pt);
			\draw [fill=black] (4.,2.) circle (4.5pt);
			\draw [fill=black] (0.,3.) circle (4.5pt);
			\draw [fill=black] (2.,3.) circle (4.5pt);
			\draw [fill=black] (4.,3.) circle (4.5pt);
			\draw [fill=black] (0.,4.) circle (4.5pt);
			\draw [fill=black] (1.,4.) circle (4.5pt);
			\draw [fill=black] (2.,4.) circle (4.5pt);
			\draw [fill=black] (3.,4.) circle (4.5pt);
			\draw [fill=black] (4.,4.) circle (4.5pt);
			\draw [fill=black] (6.,4.) circle (4.5pt);
			\draw [fill=black] (7.,4.) circle (4.5pt);
			\draw [fill=black] (8.,4.) circle (4.5pt);
			\draw [fill=black] (9.,4.) circle (4.5pt);
			\draw [fill=black] (10.,4.) circle (4.5pt);
			\draw [fill=black] (10.,3.) circle (4.5pt);
			\draw [fill=black] (9.,3.) circle (4.5pt);
			\draw [fill=black] (8.,3.) circle (4.5pt);
			\draw [fill=black] (7.,3.) circle (4.5pt);
			\draw [fill=black] (6.,3.) circle (4.5pt);
			\draw [fill=black] (6.,2.) circle (4.5pt);
			\draw [fill=black] (7.,2.) circle (4.5pt);
			\draw [fill=black] (8.,2.) circle (4.5pt);
			\draw [fill=black] (9.,2.) circle (4.5pt);
			\draw [fill=black] (10.,2.) circle (4.5pt);
			\draw [fill=black] (10.,1.) circle (4.5pt);
			\draw [fill=black] (9.,1.) circle (4.5pt);
			\draw [fill=black] (8.,1.) circle (4.5pt);
			\draw [fill=black] (7.,1.) circle (4.5pt);
			\draw [fill=black] (6.,1.) circle (4.5pt);
			\draw [fill=black] (6.,0.) circle (4.5pt);
			\draw [fill=black] (7.,0.) circle (4.5pt);
			\draw [fill=black] (8.,0.) circle (4.5pt);
			\draw [fill=black] (9.,0.) circle (4.5pt);
			\draw [fill=black] (10.,0.) circle (4.5pt);
			\draw [fill=black] (3.,3.) circle (4.5pt);
			\draw [fill=black] (1.,3.) circle (4.5pt);
		\end{scriptsize}
	\end{tikzpicture}
	\caption{The $n^2$ number of pebbles and their division}
	\label{fig:k13}
\end{figure}

\begin{fa} Using methods similar to the above, find a closed form for the sums below
	\begin{enumerate}[(a)]
		\item $2+4+\ldots+2n,$
		\item $1+2+3+\ldots+n-1+n+n-1+n-2+\ldots+2+1,$
		\item $1^2+2^2+3^2+\ldots+n^2,$
		\item $1^3+2^3+3^3+\ldots+n^3,$
		\item $f_1^2+f_2^2+\ldots+f_n^2,$ where $f_n$ denotes the $n$th Fibonacci number: $f_1=f_2=1$, $f_{n}=f_{n-1}+f_{n-2},$ $n>2$.
	\end{enumerate}
\end{fa}

\begin{fa} Let us prove and illustrate that
	\begin{enumerate}[(a)]
	\item $2+4+\ldots+2n=n(n+1),$
	\item $1+2+\ldots+n+n-1+n-2+\ldots+2+1=n^2,$
	\item $1^2+2^2+3^2+\ldots+n^2=\frac{n(n+1)(2n+1)}{6},$
	\item $1^3+2^3+3^3+\ldots+n^3=\left(\frac{n(n+1)}{2}\right)^2,$
	\item $f_1^2+f_2^2+\ldots+f_n^2=f_nf_{n+1},$ where $f_n$ denotes the $n$th Fibonacci.
	\end{enumerate}
\end{fa}

\begin{fa}  Let's examine the following equalities.
	
\[
\begin{split}
	1+3+1
	&=
	1^2+2^2,\\
	1+3+5+3+1
	&
	=2^2+3^2,\\
	1+3+5+7+5+3+1
	&=3^2+4^2.
\end{split}
\]

\end{fa}
Using the above, let's find a possible production of the sum of adjacent square numbers.

\begin{fa}
Let's make comments about the equations below and prove them.
	
\[
\begin{split}
	1 
	&=
	1^3\\
	2+4+2 
	&= 
	2^3\\
	3+6+9+6+3 
	&=
	3^3\\
	4+8+12+16+12+8+4 
	&=
	4^3
\end{split}
\]
	
\end{fa}

\subsection{Text tasks}

Solving the so-called ,,text problems'' has been an increasingly problematic area of mathematics teaching for years, it is demonstrably worse for e.g. also the graduation tasks set in this form. The comprehension of the text itself, then the setting up of the model and the interpretation of the result, are causing more and more serious problems. (The formal use of the model, i.e. the equation, system of equations, etc. solution is less.) However, the solution of this type of tasks, as well as the equation, system of equations, etc. belonging to the text, one of the most efficient ways to write it, and to create a model in general, is simple: we write the same quantity in two different ways using the known data and the unknowns introduced in a suitable way. An example of this is shown below.

\begin{fa}
	For 1.3 kg of salt solution, 8 kg of 15\% salt solution is poured, and thus a 10\% salt solution is created. What \% was the original solution?
\end{fa}

\ans{We can help you create the model by filling in a simple table, for example in the following way.
	
	\begin{center}
		 \begin{tabular}[h!]{|c|c|c|c|}
		\hline 
		& weight (kg) & concentration (\%) & amount of solute (kg) \\ 
		\hline 
		Solution 1 & 1.3 & $ x $ & $ \frac{1.3\cdot x}{100} $ \\ 
		\hline 
		Solution 2 & 0.8 & 15 & $ \frac{0.8\cdot 15}{100} $ \\ 
		\hline 
		mix & 2.1 & 10 &  \\ 
		\hline 
\end{tabular}
\end{center} }
The value belonging to an empty cell in the table, i.e. the solute content of the mixture, can be written in two ways. On the one hand, it is a 2.1 kg 10\% solution with a solute content (measured in kg) of $\frac{2.1\cdot 10}{100},$ on the other hand, the solute content of the mixture is the same as the individual with the sum of the solute content of the components, and thus $\frac{1.3\cdot x}{100}+\frac{0.8\cdot 15}{100}.$ Since we wrote the same quantity in two ways, we get:

 $$\frac{2.1\cdot 10}{100}=\frac{1.3\cdot x}{100}+\frac{0.8\cdot 15}{100}.$$
Solving the equation ($x=7$) gives the (percentage) concentration of the first component.

\subsection{The definition of a real exponential power}

In elementary school, we introduce the definition of a positive whole exponent power. Based on this, in the ninth grade we define the concepts of zero and negative integer exponents using the permanence principle. In the tenth and eleventh grades, we use the concept of the $n$th root to form the concept of a power with a rational exponent. In connection with the introduction of the exponential function, we also need the concept of a power with a real exponent. We can do all this with the help of the two-sided approximation

The function $f:\mathbb{Q}\to\mathbb{R},$ $f(x)=2^x$ is strictly monotonically increasing. It is therefore advisable to define the function $g:\mathbb{R}\to\mathbb{R},$ $g(x)=2^x$ as an extension of the function $f$ in such a way that it is also strictly monotonically increasing. For this, we need to attribute meaning to powers with irrational exponents, such as the expression $2^{\sqrt{2}}$. Consider the following two-sided approximation of $\sqrt{2}$.

$$1<\sqrt{2}<2$$ $$1.4<\sqrt{2}<1.5$$ $$1.41<\sqrt{2}<1.42$$ $$1.414<\sqrt{2}<1.415$$	$$\vdots$$
So $$2=2^1<2^{\sqrt{2}}<2^2=4$$ $$2.63\approx2^{1.4}<2^{\sqrt{2}}<2^{1.5}\approx2.83$$ $$2.657\approx2^{1.41}<2^{\sqrt{2}}<2^{1.42}\approx2.676$$ $$2.664\approx2^{1.414}<2^{\sqrt{2}}<2^{1.415}\approx2.666$$	$$\vdots$$
These intervals have exactly one point in common, which will be the value of $2^{\sqrt{2}}$. The value of the expression $a^x (a>0, a\not=1)$ can also be defined similarly, for any irrational exponent $x$ (for $0<a<1$ it is strictly monotonically decreasing, for $a>1$ it is strictly monotone increasing).

\section{Number theory}

The method of two-sided approximation can also bring many interesting results in the field of number theory.

\subsection{The number of divisors}

\begin{fa}
	Denote by $d(n)$ the number of divisors of $n$ positive integers. Let us prove that \begin{equation*}
		d(1)+d(2)+d(3)+\ldots+d(n)=n+\left[\frac{n}{2}\right]+\left[\frac{n}{3} \right]+\ldots+\left[\frac{n}{n}\right],
	\end{equation*} where $\left[a\right]$ denotes the (lower) integer part of the real number $a$, i.e. the largest of the integers not larger than $a$.
\end{fa} 

\ans{Let $n$ be a fixed positive integer, and consider a $n\times n$ table. Color the $j$-th element of the $i$-th row in the table gray if $i | j$, otherwise paint it white (Figure \ref{fig:osz1}).
	
	\begin{figure}[h!]
		\centering
		\includegraphics[width=0.3\linewidth]{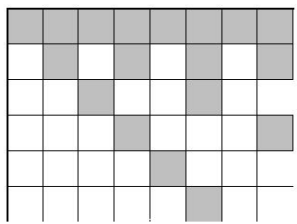}
		\vskip -3cm
		\caption{The number of divisors}
		\label{fig:osz1}
	\end{figure}
	
	Let's examine how many gray fields we have per column and per row. It is clear that in the $j$th column there is exactly $d(j)$ number of gray fields, so the total of the table is $d(1)+d(2)+d(3)+\ldots+d(n)$ contains a gray field. Counting row by row, the first row is completely gray, it has $n$ gray fields, the second row has $\left[\frac{n}{2}\right]$, the third row has $\left[\frac{n}{3} \right]$ and so on, in the last $n$th row only $1=\left[\frac{n}{n}\right]$. Thus, in the entire table $n+\left[\frac{n}{2}\right]+\left[\frac{n}{3}\right]+\ldots+\left[\frac{n}{n}\right]$ is a gray field. The two sums must match, since they both show the number of gray fields in the table, so we got: $$d(1)+d(2)+d(3)+\ldots+d(n)=n+\left[\frac{n}{2}\right]+\left[\frac{n}{3}\right]+\ldots+\left[\frac{n}{n}\right].$$}

\begin{fa}Using the above, let's estimate the $$\frac{d(1)+d(2)+d(3)+\ldots+d(n)}{n}$$ quotient.
\end{fa}

\ans{It is clear from the definition of the integralpart that if $1\leq k\leq n$ then $$\frac{n}{k}-1<\left[\frac{n}{k}\right]\leq\frac{n} {k}.$$ Thus, using our result above, we get:
	
$$(n-1)+\left(\frac{n}{2}-1\right)+\ldots+\left(\frac{n}{n}-1\right)<d(1)+d(2)+\ldots+d(n)\leq n+\frac{n}{2}+\frac{n}{3}+\ldots+\frac{n}{n},$$ 
so

$$1+\frac{1}{2}+\ldots+\frac{1}{n}-1<\frac{d(1)+d(2)+\ldots+d(n)}{n}\leq 1+\frac{1}{2}+\frac{1}{3}+\ldots+\frac{1}{n}.$$}

\sna{The interesting thing about the estimate is that both the lower and upper estimates increase beyond all limits, but their difference is 1.}

\section{Combinatorics}

Combinatorics is one of the ,,mother areas'' of two-sided addition/counting. Consider the following task.

\subsection{Binomial coefficients and what lies behind them}

\begin{fa}
	Let us prove that 
		$$1+2+\ldots+n=\binom{n+1}{2}.$$
\end{fa}

The above and similar tasks can be solved, for example, by induction, or by formal transformations using the definition and properties of binomial coefficients. None of these solution methods provide answers to the previously mentioned \textit{,,Why''}-type questions, leaving the underlying, real mathematical content unanswered.

\sna{Depending on the structure, ${n}\choose{k}$ can be defined as the number of $k$ subsets of a set with $n$ elements, or in the form ${n}\choose{k}$=$\frac{n!}{k!(n-k)!}$. In the first case, it is precisely the latter production that is displayed as an item, while in the case of the latter definition, the number of subsets of a $n$ element set with $k$ elements is obtained. The following can be used in both configurations with minor modifications.} In the following, we will solve the above task using our method.

\ans{On the right side is the number of two-element subsets of a set with $n+1$ elements. So it is enough to see that it is also on the left side. Denote the elements of the set $a_1, a_2, \ldots a_n, a_{n+1}.$ List all two-element subsets as follows. First, take the ones with the element $a_1$: $\{a_1, a_2\}$, $\{a_1, a_3\}$, $\ldots$ $\{a_1, a_n\}$, $\{a_1, a_ {n+1}\}$. Their number is $n$. Now let's take those in which the element with the smallest index is $a_2$. They are: $\{a_2, a_3\}$, $\{a_2, a_4\}$, $\ldots$, $\{a_2, a_{n+1}\}$. The number of these subsets is exactly $n-1$. Continuing this, the subsets that contain the smallest index element $a_k$: $\{a_k, a_{k+1}\}$, $\{a_k, a_{k+2}\}$, $\ldots$ $ \{a_k, a_{n+1}\}$, and their number is $n+1-k$. Thus, we actually listed all subsets, and the number of sets in each group is $n, n-1, \ldots, 2, 1$.
	
\begin{table}[h!]
		\begin{center}
			\begin{tabular}{cccccc||cc}
				
				$\{a_1, a_2\}$ & $\{a_1, a_3\}$ & $\{a_1, a_4\}$ & $\{a_1, a_5\}$ & $\ldots$ & $\{a_1, a_{n+1}\}$ & $n$ \\
				$\{a_2, a_3\}$ & $\{a_2, a_4\}$ & $\{a_2, a_5\}$ & $\ldots$ & $\{a_2, a_{n+1}\}$ & & $n-1$\\
				$\{a_3, a_4\}$ & $\{a_3, a_5\}$ & $\ldots$ & $\{a_3, a_{n+1}\}$ & & & $n-2$\\
				$\vdots$ & & & & & &  $\vdots$\\
				$\{a_{n-1}, a_n\}$ & $\{a_{n-1}, a_{n+1}\}$& & & & & 2\\
				$\{a_n, a_{n+1}\}$ & & & & & &1\\
			\end{tabular}
		\end{center}
\end{table}}

Regardless of the construction, the point where non-repeating combinations are connected to a special type of repeating permutations is of particular importance. In the following, we show an example of this, in a way that helps the discovery, using the two-sided counting method. Imagine a table that is $n-k+1$ $(0\leq k\leq n)$ steps long and $k+1$ steps wide. \begin{table}[h!]
	\begin{center}
		\begin{tabular}{|p{2mm}|p{2mm}|p{2mm}|p{2mm}|p{2mm}|p{2mm}|p{2mm}|p{2mm}|}
			\hline
			&  &  &  &  &  & C\\
			\hline
			&  &  &  &  &  & \\
			\hline
			&  &  &  &  &  & \\
			\hline
			&  &  &  &  &  & \\
			\hline
			M &  &  &  &  &  & \\
			\hline
		\end{tabular}
		\caption{The mouse and the cheese}
	\end{center}
\end{table} 

There is currently a mouse sitting in the lower left corner, and a piece of cheese in the upper right corner of the board. Our mouse can only move up or to the right, always to an adjacent field. How many ways can our mouse get to the cheese?\\ \indent First, let's see how many steps our mouse needs. Since it is currently in the first space of the first row, it needs to move $n-k$ more to the right and move up $k$ spaces. Our \textit{clever} mouse can think in two ways. \\ \indent \textbf{1.} Since he needs a total of $n-k$ steps to the right and $k$ up, he creates $n$ number of cards and draws $\rightarrow$ on $n-k$ and $\uparrow$ on $k$. Now, as many ways as you can line up these arrows, you can get to the cheese in as many ways. And the number of these row arrangements is exactly $=\frac{n!}{k!\cdot(n-k)!}$, since it has $n$ number of cards, of which $k$ and $n-k$ are the same. \\ \indent \textbf{2.} After thinking about this, our mouse finds it too tiring, so it thinks like this: If I number my steps from one to $n$ and select the ones I want to take upwards, it is clear that I have to take the unselected steps to the right. Thus, there are as many ways to the cheese as there are ways in which I can choose $k$ from among the numbers $1, 2, \ldots, n$, and the number of these choices is exactly $\binom{n}{k}$. \\ \indent Thus we provied that $$\binom{n}{k}=\frac{n!}{k!\cdot(n-k)!}.$$

The difficulty of a task, its success in solving it, and its popularity among students depend not only on the problem itself (so all of the above are subjective), but also on the way it is presented. To illustrate this, consider the following task.

\begin{fa}
	 \begin{enumerate}[(a)]
	 	\item Add the elements in rows $0, 1, 2, \ldots$ of Pascal's triangle. What do we experience? Let's state our conjecture in general terms and prove it.
	 	\item Let us prove that the number of subsets of a set with $n$ elements is $2^n$.
	 	\item Let us prove that the sum of the elements in the $n$th row of Pascal's triangle is $2^n$.
	 	\item Prove that $\binom{n}{0}+\binom{n}{1}+\ldots+\binom{n}{n}=2^n$.
	 \end{enumerate}
 \end{fa}

\ans{In the row $n$ of Pascal's triangle, or on the left side of equality (d), there is the number of all subsets of a set with $n$ elements. Let the elements of the set be $a_1, a_2, \ldots a_n$. For each subset, it is true whether $a_1$ is an element or not. Thus, we can decide in two ways regarding the element $a_1$ (either we include it in the given subset or not). Whatever we decide, regardless of this decision, we have to decide again whether to include it in the subset in the case of element $a_2$. Thus, in the case of the first two elements, there are already $2\cdot 2=4$ possible choices. The procedure can be continued, and no matter how we decided in the case of the first $k$ element, the number of possibilities is doubled if the next element is taken into account. Thus, when we reach the $n$th element, we will have $2^n$ different subsets.
	
The above method becomes clearer and more understandable if we list the elements of the set and write a $+$ or $-$ sign (possibly 0 or 1) under each element, depending on whether it is included in the given subset, or not.
	
\begin{table}[h!]
		\begin{center}
			\begin{tabular}{cccccc}
				
				$a_1,$ & $a_2,$ & $a_3,$ & $\ldots$ & $a_{n-1},$ & $a_n$  \\
				$+$ & $-$ & $-$ & & $+$ & $-$\\
			\end{tabular}
		\end{center}
\end{table}
	
Thus, there can be two types of signs under each element independently, so the number of possibilities is: $2\cdot 2\cdot \ldots \cdot 2=2^n$}.

\begin{fa} Let us prove and illustrate that
	\begin{enumerate}[(a)]
		\item $\binom{n+1}{k}=\binom{n}{k}+\binom{n}{k-1},$
		\item if $2\leq k\leq n$, then $\binom{n+1}{k}=\binom{n-1}{k-2}+2\cdot\binom{n-1}{k-1}+\binom{n-1}{k},$
		\item if $0\leq j\leq k\leq n$, then $\binom{n+1}{k}=\sum_{i=0}^j\binom{j}{i}\binom{n+1-j}{k-i}.$
	\end{enumerate}
\end{fa}

\ans{\begin{enumerate}[(a)]
		\item Let's imagine that we have to choose a delegation of $k$ members in a class with $n+1$ members. Of course, we can do this in a $\binom{n+1}{k}$ way. But if we pay special attention to the fate of Ábel (who is, of course, a member of the class), we can find another solution to the same problem. Our situation is not much more complicated now, Ábel is either included in the delegation or not. If so, then he is already a member, so $k-1$ more delegates must be chosen from the remaining $n$ members of the class, and this can be done in the $\binom{n}{k-1}$ way. And if Ábel is not a member of the delegation, then all $k$ delegates must be selected from the other $n$ members of the class, and we have $\binom{n}{k}$ options for this. In this way, delegates can be chosen in total, according to a different point of view, in a $\binom{n}{k}+\binom{n}{k-1}$ manner.
		\item Now let's select separately the cases where, of the two members of the class, Ábel and Balázs, both exactly one or neither of them is a member of the committee.
		\item Let us now monitor the fate of student $j$ in the class.
\end{enumerate}}

The two-sided counting method not only shows the essence, but also provides simple, beautiful solutions that do not involve calculations or formal transformations, sometimes with surprising interpretations. With this in mind, let's solve the following task.

\begin{fa}\begin{enumerate}[(a)]
	\item $2^n\cdot\binom{n}{0}+2^{n-1}\cdot\binom{n}{1}+
	2^{n-2}\cdot\binom{n}{2}+\ldots+2\cdot\binom{n}{n-1}+1\cdot\binom{n}{n}=3^n,$
	\item $\sum_{k=0}^n\sum_{m=0}^k\binom{n}{k}\binom{k}{m}=3^n.$
		\end{enumerate}
\end{fa}

The relationship between Pascal's triangle containing binomial coefficients and Fibonacci numbers can be discussed using our method in the following way. Examining Pascal's triangle a little more closely, we can notice that by summing up the lines of Pascal's triangle in a suitable way, we get the correct Fibonacci numbers.

\begin{center}
	\includegraphics[width=0.4\linewidth]{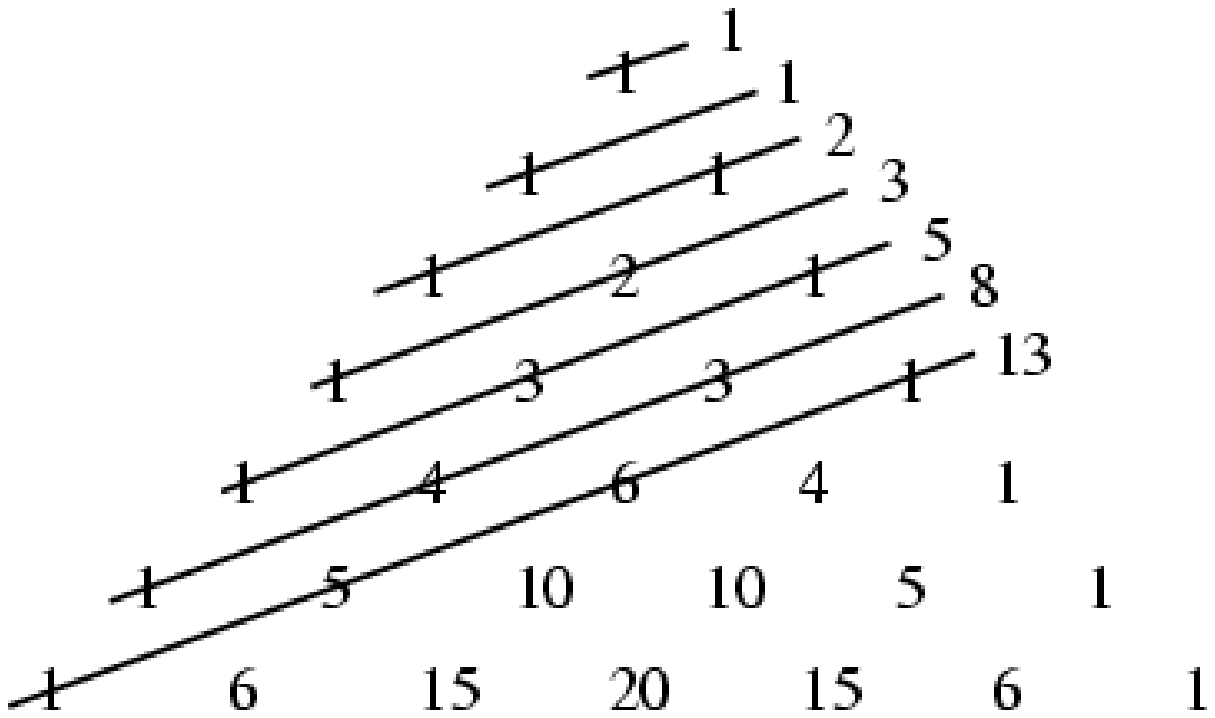}
\end{center}

The observation itself can be surprising, and even put into a formula, we get the following, not very friendly task.

\begin{fa}

Let us prove that $$f_n=\sum_{k=0}^{\left[{{n-1} \over 2}\right]} \binom{n-k-1}{k},$$ where $f_n$ denotes the $n$th Fibonacci number.
	
\end{fa}

Neither the method of determination, nor the possible solution methods provide an answer to the reason for this surprising correlation. Starting from a little further away, however, we can find a solution that answers all questions. So let's solve the following exercise (in two ways).
	
\begin{fa} 
	
We want to paint all floors of a house blue or white, but two adjacent floors cannot be blue at the same time. How many different colors are possible for a house with $n$ floors?
	
\end{fa}

\ans{\textbf{(1.)} So let $h_n$ denote the number of possible colorings of a $n$ storey house. In the case of a one-floor house, we can color it in two ways, or we can paint the single floor white or blue. In the case of two floors, we have three options for coloring, as follows: white-white, white-blue, blue-white. In the case of a house with $n$ floors, the $n$th floor is either white or blue. If it is white, then the floor $n-1$ below it can be colored arbitrarily, i.e. the number of properly colored houses with $n$ floors is the same as the number of properly colored houses with $n-1$ floors $(h_{n-1})$. If the color of the last $n$th floor is blue, then the $n-1$th floor is necessarily white, so the remaining $n-2$ floors must be colored, which is possible in exactly as many ways as there are $n-2$ floors house can be painted according to the conditions $(h_{n-2})$. Then: $h_1=2$, $h_2=3$, $h_n=h_{n-1}+h_{n-2}$ if $n \geq 3.$ That is: $h_n=f_{n+2} $.
	
\medskip 	

\textbf{(2.)} Let the $k$ floors of the $n$ floors house be blue. Then there is a total of $n-k+1$ space between, below and above the remaining $n-k$ white floors, we can place our $k$ blue floors here, which is $\binom{n-k+1}{k}$ possible. The number of blue floors can vary from 0 until the above placement is possible, i.e.: 
\[\begin{split}
k 
&\leq 
n-k+1\\
2k
&\leq
n+1\\
k 
&\leq
\left[{{n +1} \over 2}\right]
\end{split}\]  
So we got:

$$h_n=\sum_{k=0}^{\left[{{n+1} \over 2}\right]} \binom{n-k+1}{k}.$$ 
Rewriting the relation from $n+2$ to $n$, we are ready.}

\begin{fa}
	Let us prove and illustrate that
	\begin{enumerate}[(a)]
		\item $n\binom{n-1}{k}=k\binom{n}{k},$
		\item $\binom{n}{l}\binom{n-l}{k-l}=\binom{k}{l}\binom{n}{k},$
		\item $\binom{n}{k}+\binom{n-1}{k}+\ldots+\binom{k}{k}=\binom{n+1}{k+1}.$
	\end{enumerate}
\end{fa}

\begin{fa}
	How many head-write sequences of length $n$ are there in which two heads cannot be next to each other?
\end{fa} 

\subsection{Partitions}

By partitions of a positive integer $n$, we mean writing the number $n$ as a sum of positive integers, where the order of the terms to be added does not matter.

\begin{fa}
	
Let us prove that if $n, k \in \mathbb{N}^+$ then the number of partitions of $n$ in which each member is at most $k$ is the same as the number of partitions of $n$ containing at most $k$ members.

\end{fa}

\ans{The number of partitions of the number $n$ can be illustrated naturally as follows. Place unit squares next to or above each other so that each partition consists of a number of unit squares corresponding to the size of the partition.
	
\begin{center}
	\includegraphics[width=0.2\linewidth]{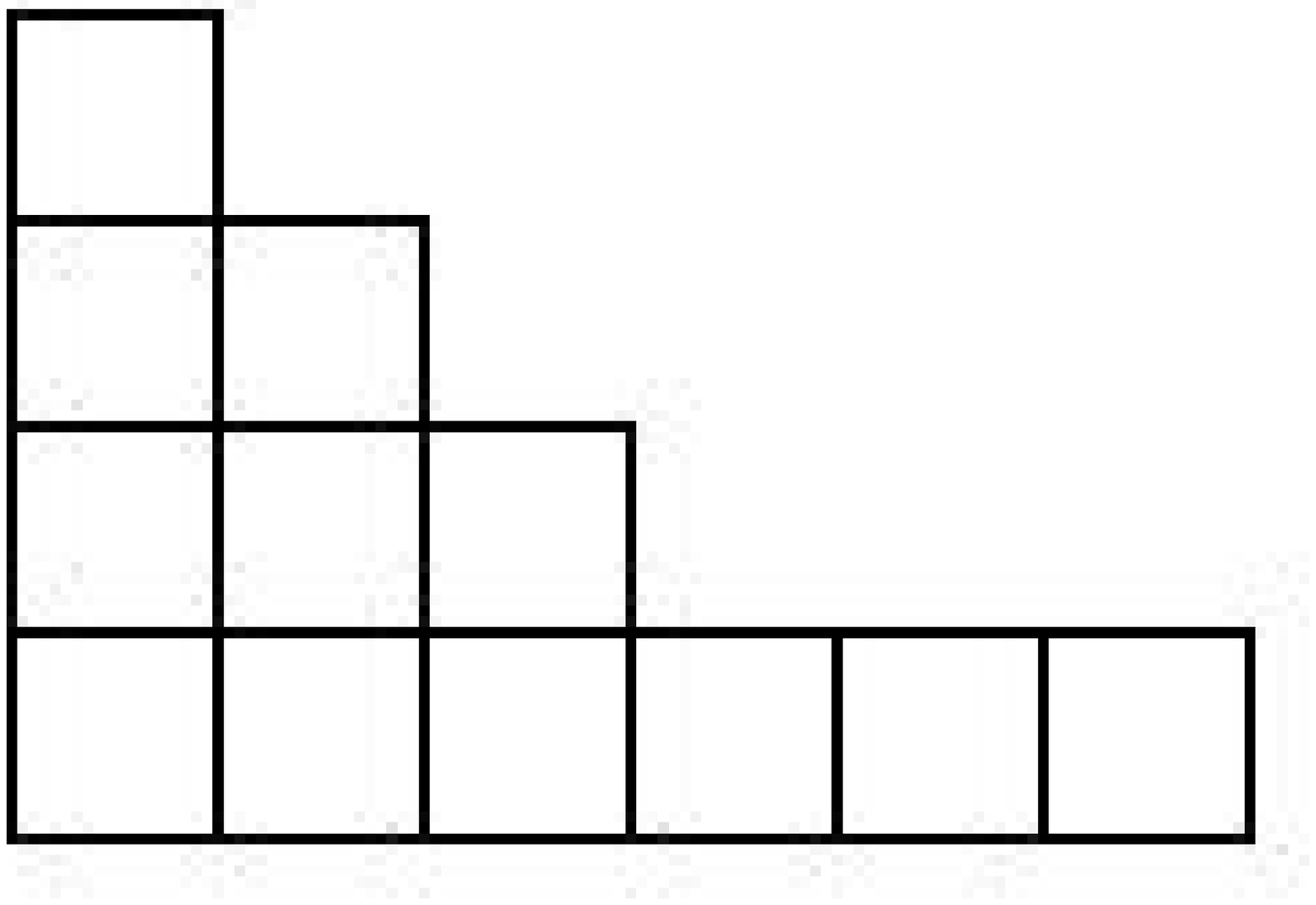}
\end{center}
\vskip -2.5cm
	
Let's rotate the figure by $90^{\textdegree}$! (Or simply look at the same figure from a different direction.)
	
\begin{center}
	\includegraphics[width=0.2\linewidth]{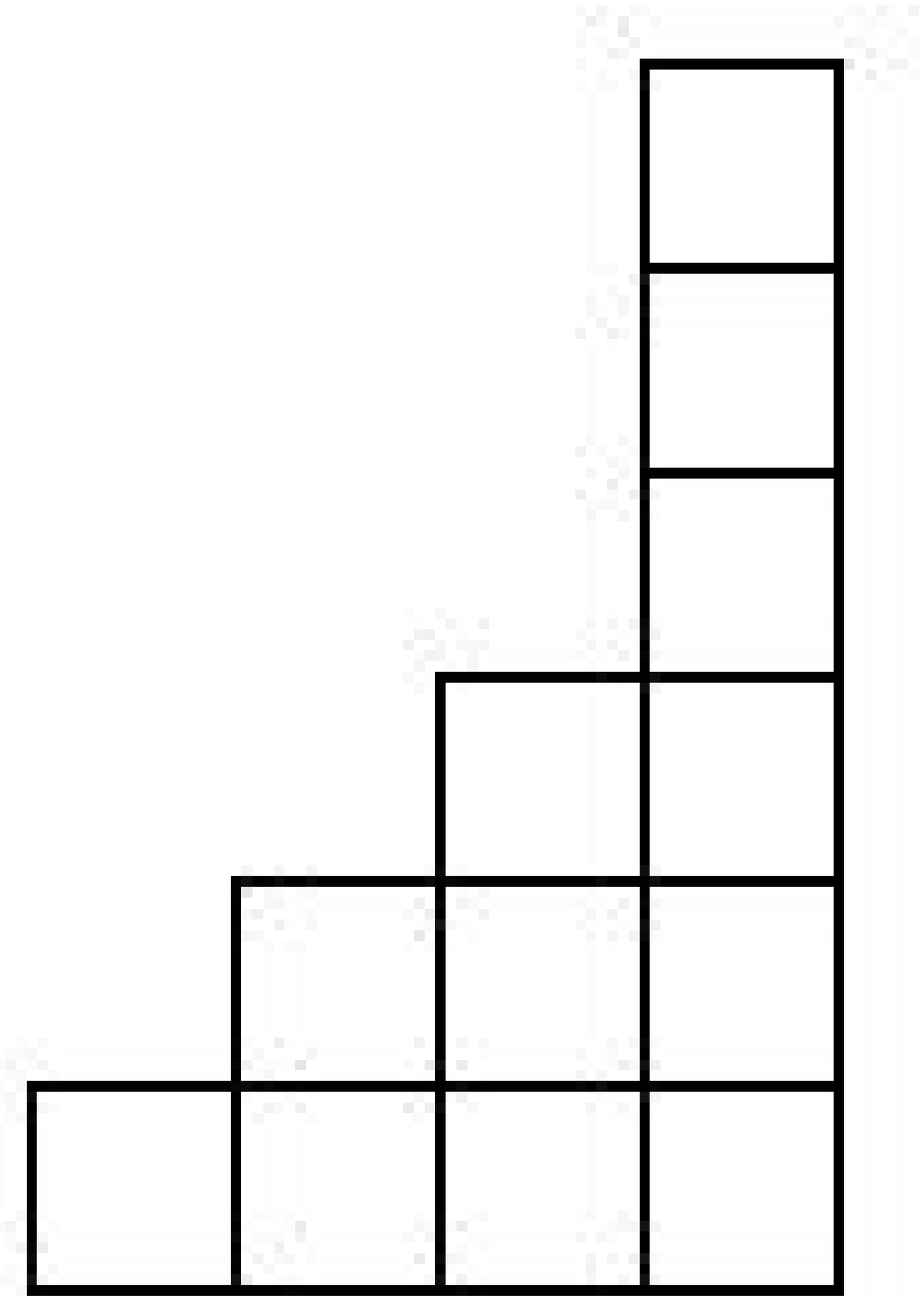}
\end{center}
	
\vskip -1cm

If all members of the partition in the first figure are at most $k$, this means that the rotated figure shows a partition with at most $k$ members. Therefore, every diagram belonging to a partition containing at most $k$ is also a diagram of exactly one partition with at most $k$ members. With this observation, we realized our claim.}

\section{Analysis}

In the following, we show some examples of the appearance and application possibilities of our method from the field of analysis.

\subsection{Police (or sandwich) principle} 

One of the relations of elementary analysis concerning the limit value of series is the theorem also known as the police (sandwich) principle:

\begin{theorem}
Let $a_n, b_n, c_n$ be sequences of real numbers for which (almost) for every $n\in \mathbb{N}$ $a_n\leq b_n\leq c_n$. Then if $\lim_{n\to\infty}a_n=\lim_{n\to\infty}c_n=a$, then the sequence $b_n$ is also convergent, and $\lim_{n\to\infty}b_n= a$.
\end{theorem}

\begin{fa}
Determine the limit of the sequence $b_n=\sqrt[n]{3^n+5^n}$.
\end{fa} 

\ans{Consider the following estimate, which holds for all $n\in\mathbb{N}$: $$a_n=\sqrt[n]{5^n}\leq\sqrt[n]{3^n+5^n} \leq\sqrt[n]{5^n+5^n}=c_n.$$} Since $a_n=5$ and $c_n=\sqrt[n]{2}\cdot5$, and so $\lim_{ n\to\infty}a_n=5,$ and $\lim_{n\to\infty}c_n=5,$ therefore $$\lim_{n\to\infty}b_n=\lim_{n\to\infty} \sqrt[n]{3^n+5^n}=5.$$

\subsection{Infinite series - A geometric series}

With the help of our method, (certain types) of infinite sums become manageable and negotiable with the help of elementary tools.

\begin{fa}
	In the past, they sold chocolate that had a slip in its paper, and for ten slips you could get another bar of chocolate. How many bars of edible chocolate is one such bar actually worth? (The task comes from László Kalmár (1905-1976).)
\end{fa}

\ans{
Each bar of chocolate has a coupon, which according to the rules is worth one tenth of chocolate. But this tithe chocolate also contains a tithe label, which is also worth all the chocolates, and so on. So a bar of chocolate is actually
	
$$1+\frac{1}{10}+\frac{1}{10^2}+\frac{1}{10^3}+\ldots+\frac{1}{10^n }+\ldots$$
worth of chocolate.

Let's look at our chocolate from a different perspective. Pisti goes into the store with 9 coupons and asks the saleswoman, with whom she is on good terms, since she often buys chocolate from her, to give her a chocolate from the shelf for a moment, and she will pay for it without waiting. Pisti, as soon as he received the chocolate, opens it and takes out the missing tenth coupon for payment, and then pays for his chocolate with the ten coupons. So he got 10 for the price of 9 chocolates, that is, one chocolate is actually worth $\frac{10}{9}$ chocolates. With this we realized:

$$1+\frac{1}{10}+\frac{1}{10^2}+\frac{1}{10^3}+\ldots+\frac{1}{10^n}+\ldots=\sum_{n=0}^\infty\frac{1}{10^n}=\frac{10}{9}.$$
}

\subsection{A non-geometric series - a geometric model}

\begin{fa}
	
Determine the sum of the row\footnote{Richard Swineshead, an English mathematician, dealt with this line for the first time during the investigation of a physical problem in the XIV. century.} \[\sum_{n=0}^\infty\frac{n}{2^n}.\]

\end{fa}

\ans{Consider the following figure (Figure \ref{fig:njsor}). 
	
\begin{figure}[h!]
		\centering
		\includegraphics[width=0.55\linewidth]{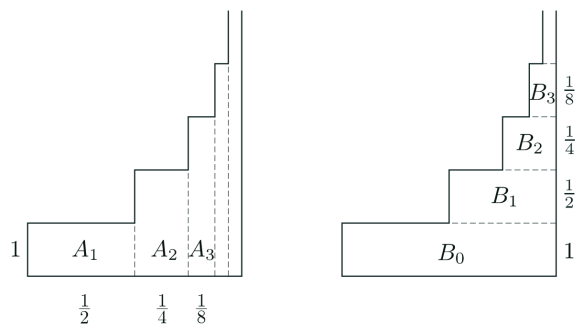}
		\caption{The amount is two ways}
		\label{fig:njsor}
\end{figure}
		
The sum of the areas of the rectangles shown on the left: $\frac{1}{2}+\frac{2}{2^2}+\frac{3}{2^3}+\ldots+\frac{n}{2^n}+\ldots,$ it is clear that it is equal to the sum of the areas of the rectangles on the right - the two figures can be shifted into each other - which: $1+\frac{1}{2}+\frac{1}{2^2}+\frac{1}{2^3}+\ldots+\frac{1}{2^n}+\ldots.$ And this sum - as it is the sum of a convergent geometric series - can be easily calculated. This is how we got:
 
$$\sum_{n=0}^\infty\frac{n}{2^n}=1+\frac{1}{2}+\frac{1}{2^2}+\frac{1}{2^3}+\ldots+\frac{1}{2^n}+\ldots=2.$$}
	
\sna{The sum on the right can also be determined, even without knowledge of convergent geometric series, with the help of the following problem.}
	
\begin{fa}
	We have a circular cake that we want to eat completely so that we can eat it every day of our endlessly long lives. How should we proceed?
\end{fa}
	
\ans{On the first day, eat half of the cake, on the second day, half of the remaining part, and so on. Thus, on the $n$th day, we eat the $\frac{1}{2^n}$th part of the cake, i.e.: $\frac{1}{2}+\frac{1}{2^2} +\frac{1}{2^n}+\ldots+\ldots=1.$} The series in the task can be summed up in a different way, in a way that does not exceed advanced high school knowledge, the sum can be calculated as the sum of the sum of the next, imaginatively written infinite number of geometric series.

\[
\begin{split}
\dfrac{1}{2}+\dfrac{1}{2^2}+\dfrac{1}{2^3}+\dfrac{1}{2^4}
&+
\ldots\\
\dfrac{1}{2^2}+\dfrac{1}{2^3}+\dfrac{1}{2^4}
&+
\ldots\\
\dfrac{1}{2^3}+\dfrac{1}{2^4}
&+
\ldots\\
&\vdots
\end{split}
\]
Summing up the sums of the above geometric series, $1+\frac{1}{2}+\frac{1}{2^2}+\frac{1}{2^3}+\ldots+\frac{1}{2^n}$ geometric series whose sum is 2.

\subsection{The definite integral and the area under the graph}

Integral calculus is one of the areas of analysis that is also important from the point of view of applications. Already in the famous Moscow Scroll created around 1800 BC, the basic idea of the definite integral calculus can be found in connection with the calculation of the volumes of truncated cones and truncated pyramid. A more advanced version of the method can be found among the ancient Greeks.

The so-called infinitesimals (infinitely small quantities) were first used by Archimedes, with which he achieved significant results. However, it should be mentioned that
even Archimedes himself did not consider his own proofs to be accurate. Integral calculus in today's sense was developed following the work of Newton (1642-1727) and Leibniz, as well as Cauchy and Riemann (1826-1866).

The calculation of the ,,size'' of various plane domes is of serious practical importance. With polygons, this is usually not a problem if you have the right data, but determining (verifying) the area of a circle is not necessarily easy either. Next, we describe a procedure developed by the aforementioned Archimedes to determine the area of a parabolic slice. This is commonly called the method of two-sided approximation.

\begin{fa}
	Calculate the area under the graph of the parabola with the equation $y=x^2$ on the closed interval $[0,1]$.
\end{fa}

\begin{meg*}
		
Intuitively, we can immediately give an estimate of the area, since it is at least zero, of course larger than it, or less than half the area of the unit square. However, we would like to define the area more precisely than that. We can do this as follows.
	
\vspace{0.2cm}

\begin{figure}[H]
			\centering
			\includegraphics[width=0.4\linewidth]{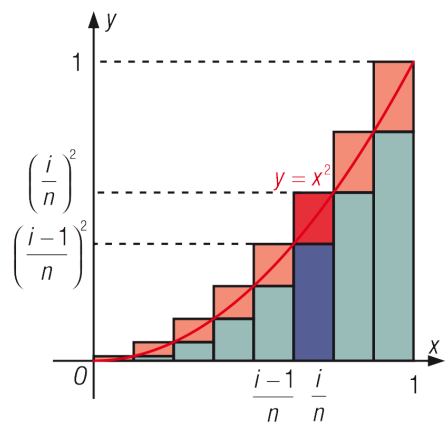}
			\label{fig:integral}
			\caption{Two-sided approximation of the area under the graph}
\end{figure}

We denote the area to be searched by $T$, which we already know is $0<T<1/2$. Divide the interval $[0,1]$ into $n$ equal parts, where $n\in\N$. The area under the curve can be estimated from the bottom and top of the $i$th subinterval using the area of a rectangle. The area of the rectangles belonging to the $i$th subinterval is (based on the Figure \ref{fig:integral}):

	\[
	\dfrac{1}{n}\cdot \left(\dfrac{i-1}{n}\right)^2=\dfrac{(i-1)^2}{n^3}\text{, and } \dfrac{1}{n}\cdot \left(\dfrac{i}{n}\right)^2=\dfrac{i^2}{n^3},
	\]
thus, the area $T$ under the graph on the interval $[0,1]$ is estimated simultaneously from the bottom and top:
\[
	\dfrac{1}{n^3}(1^2+2^2+\ldots+(n-1)^2)\leq T\leq \dfrac{1}{n^3}(1^2+2^2+\ldots+n^2).
\]
Then it refers to the closed form of the square of the first positive integer $n$
	
	\[
	1^2+2^2+\ldots+n^2=\dfrac{n(n+1)(2n+1)}{6}.
	\]
with the help of a relation, we can give the following estimate:
	
	\[
	\begin{split}
		\dfrac{(n-1)n(2n-1)}{6n^3}
		& \leq T \leq \dfrac{n(n+1)(2n+1)}{6n^3},\\
		\left(1-\dfrac{1}{n}\right)\left(\dfrac{1}{3}-\dfrac{1}{6n}\right)
		& \leq T \leq \left(1+\dfrac{1}{n}\right)\left(\dfrac{1}{3}+\dfrac{1}{6n}\right).
	\end{split}
	\]
	Then using the limit value
	
	\[\lim_{n\to \infty}\dfrac{1}{n}=0,\]
it can be seen that the limit value of the series on both sides of the inequality is $\dfrac{1}{3}$, that is, $T=\dfrac{1}{3}$.
	
\end{meg*}

\begin{fa}
	Given the function $f:[0,a]\rightarrow\R, f(x)=x^3$, where $a>0$. Determine the area of the closed plane section under the graph of the function.
\end{fa}

The success of determining the area in this way also depended on the calculation of the lower and upper approximate amounts. 

\begin{def*}
	
A bounded function $f$ interpreted on a closed interval $[a,b]$ is called integrable exactly if there is only one number that is neither smaller than any lower approximate sum of the function $f$ nor larger than any upper approximate sum. This number is called the definite (Riemann) integral of the function $f$ on $[a,b]$, and

\[
\int_{a}^{b} f(x)\dd x
\]
marked as functions with this property are usually called Riemann-integrable functions.
\end{def*} 
Using this notation based on the above
\[
\int_{0}^{1} x^2 \dd x=\dfrac{1}{3}
\]
occurs.

\section{Geometry}

\subsection{Analytic geometry - Dot product of vectors}

One of the simplest, textbook ways of calculating the inclination angle of vectors given with their coordinates in the Cartesian coordinate system results from writing the dot product of the vectors in two ways. Let the vectors $\overrightarrow{a}(a_1; a_2)$ and $\overrightarrow{b}(b_1; b_2)$ be given on the coordinate plane, denote the inclination angle of these vectors by $\alpha$. Then the dot product of the two vectors is by definition: $$\overrightarrow{a}\cdot\overrightarrow{b}=|\overrightarrow{a}|\cdot|\overrightarrow{b}|\cdot\cos\alpha=\sqrt{ a_1^2+a_2^2}\cdot\sqrt{b_1^2+b_2^2}\cdot\cos\alpha.$$ Expressed with coordinates,
$$\overrightarrow{a}\cdot\overrightarrow{b}=a_1\cdot b_1+a_2\cdot b_2.$$
It follows from the two-sided approach of the dot product that

\begin{equation}\label{ska1}
\sqrt{a_1^2+a_2^2}\cdot\sqrt{b_1^2+b_2^2}\cdot\cos\alpha=a_1\cdot b_1+a_2\cdot b_2.
\end{equation}
Hence, expressing the cosine of the inclination angle of the two vectors:
$$\cos\alpha=\frac{a_1\cdot b_1+a_2\cdot b_2}{\sqrt{a_1^2+a_2^2}\cdot\sqrt{b_1^2+b_2^2}}.$$

\sna{From the equation (\ref{ska1}), taking into account that $\cos\alpha\leq 1$, the Cauchy-Schwarz-Bunyakovskii inequality is easily obtained (in two terms): Arbitrary $a_1; a_2; b_1; b_2$ for real numbers: $$\sqrt{a_1^2+a_2^2}\cdot\sqrt{b_1^2+b_2^2}\geq a_1\cdot b_1+a_2\cdot b_2 .$$ }

\subsection{Lattice geometry - Pick's formula}

One of the interesting areas of geometry is lattice geometry, which provides the opportunity to approach the geometry of number theory theorems. Points of the Cartesian coordinate system with integer coordinates are called grid points. A lattice polygon is a polygon whose vertices are all lattice points. (It is interesting, for example, that there are only equilateral lattice polygons with an even number of sides.) A lattice triangle is called empty if, apart from its vertices, it contains no lattice points either on its border or in its interior. A central relation of lattice geometry is the Pick formula, which defines the area of grid polygons using the number of ,,connected'' grid points. We achieve this using our method below. We use the fact that the area of every empty lattice triangle is $\frac{1}{2}$. Solve the following problem first.

\begin{fa}
	Prove that if there are $h$ grid points on the border of a grid polygon and $b$ grid points inside, then the polygon can be divided into $h+2b-2$ empty grid triangles.
\end{fa}

\ans{The division of grid polygons into empty grid triangles is not clear, but the number of empty grid triangles required for the resolution is independent of the method of resolution, and the sought relationship depends on this surprising statement, the proof of which rests on our method. Let the number of empty grid triangles generated during a given resolution (Figure \ref{fig:rg2}) be $n$. Calculate the sum of the interior angles of the triangles included in the resolution in two ways.\\ \textbf{(1.)} The first ,,mode'' is very simple: the sum of the interior angles of a $n$ triangle is $n\cdot 180$\textdegree.\\ \textbf{(2.)} On the other hand, the sum is obtained from the $360$\textdegree angles formed at the vertices around the interior points, on the other hand, looking at the polygon as a $h$-angle containing 180\textdegree angles from the total $(h-2)180$\textdegree angles created at the ,,vertices''. \begin{figure}[h!]
		\centering
		\definecolor{zzttqq}{rgb}{0.6,0.2,0.}
		\definecolor{qqwuqq}{rgb}{0.,0.39215686274509803,0.}
		\definecolor{cqcqcq}{rgb}{0.7529411764705882,0.7529411764705882,0.7529411764705882}
		\begin{tikzpicture}[scale=0.8][line cap=round,line join=round,>=triangle 45,x=1.0cm,y=1.0cm]
			\draw [color=cqcqcq,, xstep=1.0cm,ystep=1.0cm] (-0.5,-0.5) grid (8.5,5.5);
			\clip(-0.5,-0.5) rectangle (8.5,5.5);
			\draw [shift={(6.,2.)},line width=2.pt,color=qqwuqq,fill=qqwuqq,fill opacity=0.10000000149011612] (0,0) -- (45.:0.5694822753146987) arc (45.:296.565051177078:0.5694822753146987) -- cycle;
			\fill[line width=2.pt,color=zzttqq,fill=zzttqq,fill opacity=0.05000000074505806] (0.,0.) -- (1.,5.) -- (5.,4.) -- (8.,4.) -- (7.,3.) -- (6.,2.) -- (7.,0.) -- (4.,1.) -- (3.,0.) -- (3.,3.) -- cycle;
			\draw [line width=2.pt] (0.,0.)-- (1.,5.);
			\draw [line width=2.pt] (1.,5.)-- (5.,4.);
			\draw [line width=2.pt] (5.,4.)-- (8.,4.);
			\draw [line width=2.pt] (8.,4.)-- (6.,2.);
			\draw [line width=2.pt] (6.,2.)-- (7.,0.);
			\draw [line width=2.pt] (7.,0.)-- (4.,1.);
			\draw [line width=2.pt] (4.,1.)-- (3.,0.);
			\draw [line width=2.pt] (3.,0.)-- (3.,3.);
			\draw [line width=2.pt] (3.,3.)-- (0.,0.);
			\draw [line width=2.pt,color=zzttqq] (0.,0.)-- (1.,5.);
			\draw [line width=2.pt,color=zzttqq] (1.,5.)-- (5.,4.);
			\draw [line width=2.pt,color=zzttqq] (5.,4.)-- (8.,4.);
			\draw [line width=2.pt,color=zzttqq] (8.,4.)-- (7.,3.);
			\draw [line width=2.pt,color=zzttqq] (7.,3.)-- (6.,2.);
			\draw [line width=2.pt,color=zzttqq] (6.,2.)-- (7.,0.);
			\draw [line width=2.pt,color=zzttqq] (7.,0.)-- (4.,1.);
			\draw [line width=2.pt,color=zzttqq] (4.,1.)-- (3.,0.);
			\draw [line width=2.pt,color=zzttqq] (3.,0.)-- (3.,3.);
			\draw [line width=2.pt,color=zzttqq] (3.,3.)-- (0.,0.);
			\begin{scriptsize}
				\draw [color=black] (0.,0.) circle (3.5pt);
				\draw [color=black] (1.,5.) circle (3.5pt);
				\draw [color=black] (5.,4.) circle (3.5pt);
				\draw [color=black] (8.,4.) circle (3.5pt);
				\draw [color=black] (6.,2.) circle (3.5pt);
				\draw [color=black] (7.,0.) circle (3.5pt);
				\draw [color=black] (4.,1.) circle (3.5pt);
				\draw [color=black] (3.,0.) circle (3.5pt);
				\draw [color=black] (3.,3.) circle (3.5pt);
				\draw [color=black] (1.,4.) circle (3.5pt);
				\draw [color=black] (1.,3.) circle (3.5pt);
				\draw [color=black] (1.,2.) circle (3.5pt);
				\draw [color=black] (1.,1.) circle (3.5pt);
				\draw [color=black] (2.,2.) circle (3.5pt);
				\draw [color=black] (2.,3.) circle (3.5pt);
				\draw [color=black] (2.,4.) circle (3.5pt);
				\draw [color=black] (3.,4.) circle (3.5pt);
				\draw [color=black] (4.,4.) circle (3.5pt);
				\draw [color=black] (4.,3.) circle (3.5pt);
				\draw [color=black] (4.,2.) circle (3.5pt);
				\draw [color=black] (5.,1.) circle (3.5pt);
				\draw [color=black] (5.,2.) circle (3.5pt);
				\draw [color=black] (5.,3.) circle (3.5pt);
				\draw [color=black] (6.,3.) circle (3.5pt);
				\draw [color=black] (6.,4.) circle (3.5pt);
				\draw [color=black] (7.,4.) circle (3.5pt);
				\draw [color=black] (7.,3.) circle (3.5pt);
				\draw [color=black] (6.,1.) circle (3.5pt);
			\end{scriptsize}
		\end{tikzpicture}
		\caption{Grid polygon and its resolution}
		\label{fig:rg2}
	\end{figure}
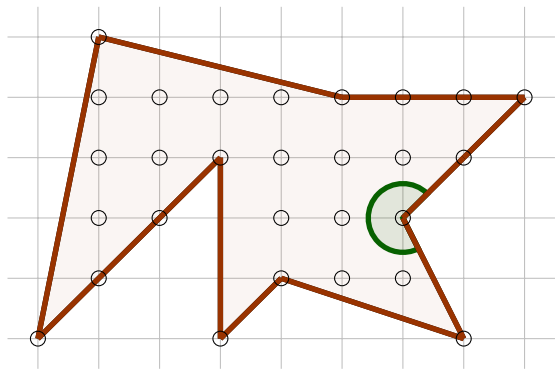

Then the joint result of the two types of counting gives $n\cdot 180=(h-2)180+b\cdot 360,$ and thus $n=h-2+2b$.} Now using the fact that the area of every empty grid triangle is: $t=\frac{1}{2}$, if there is a grid point $h$ on the boundary of a grid polygon and $b$ inside it, then the area of the polygon is: $t=\frac{h }{2}+b-1$, this is the so-called Pick formula.

\subsection{Synthetic geometry - Concept of area (Jordan measure)}

In this point, we will define the area and volume of certain subsets of the plane and the space without claiming to be complete, in a different approach from what we have seen before. The concept of area and its measurement is already introduced in the lower grades of elementary school. Then, instead of axioms, we say the following: The basic unit of the area is 1 square meter ($m^2$), the area of a square with a side length of 1 $m$. 

\newpage

\begin{fa}
	Find a method by which we can decide which of the two planar domes shown here (Figure \ref{fig:paca1}) is larger.
	
	\begin{figure}[h!]
		\centering
		\includegraphics[width=0.3\linewidth]{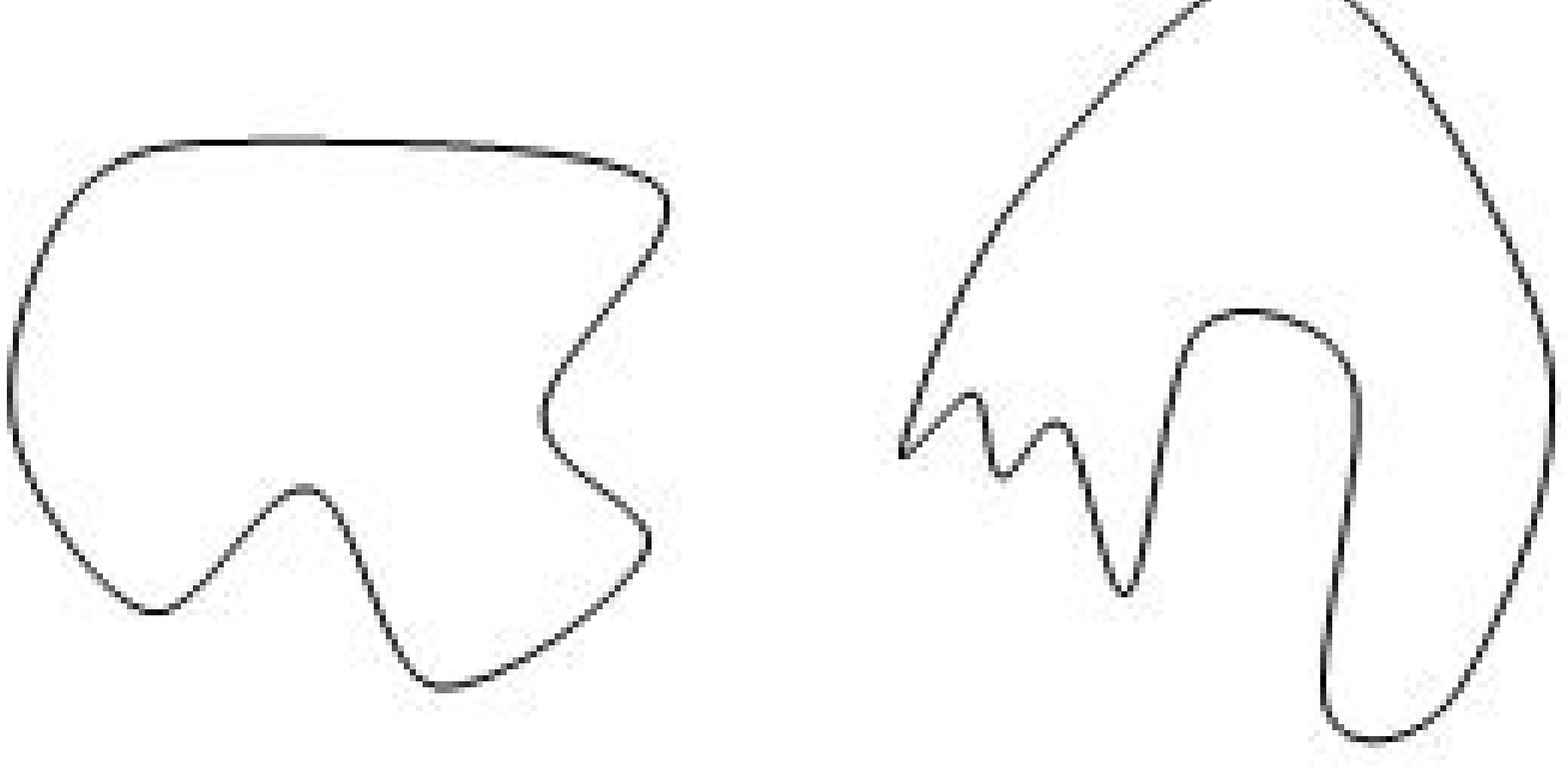}
		\vskip-4cm
		\caption{Area of planeidom}
		\label{fig:paca1}
	\end{figure}
	
\end{fa}

\ans{The solution at this age can be imagined in the following way, for example: Cut out the planar domes, then place them on a square grid with a side length of 1 cm (millimeter paper) and count how many square grids are inside each plane dome, and how many square grids can cover each plane dome juice. If this estimate (regarding the area) does not help, let's examine a square grid consisting of squares with a side length of 0.5 cm, and so on (Figure \ref{fig:p2})!
	
	\begin{figure}[h!]
		\centering
		\includegraphics[width=0.3\linewidth]{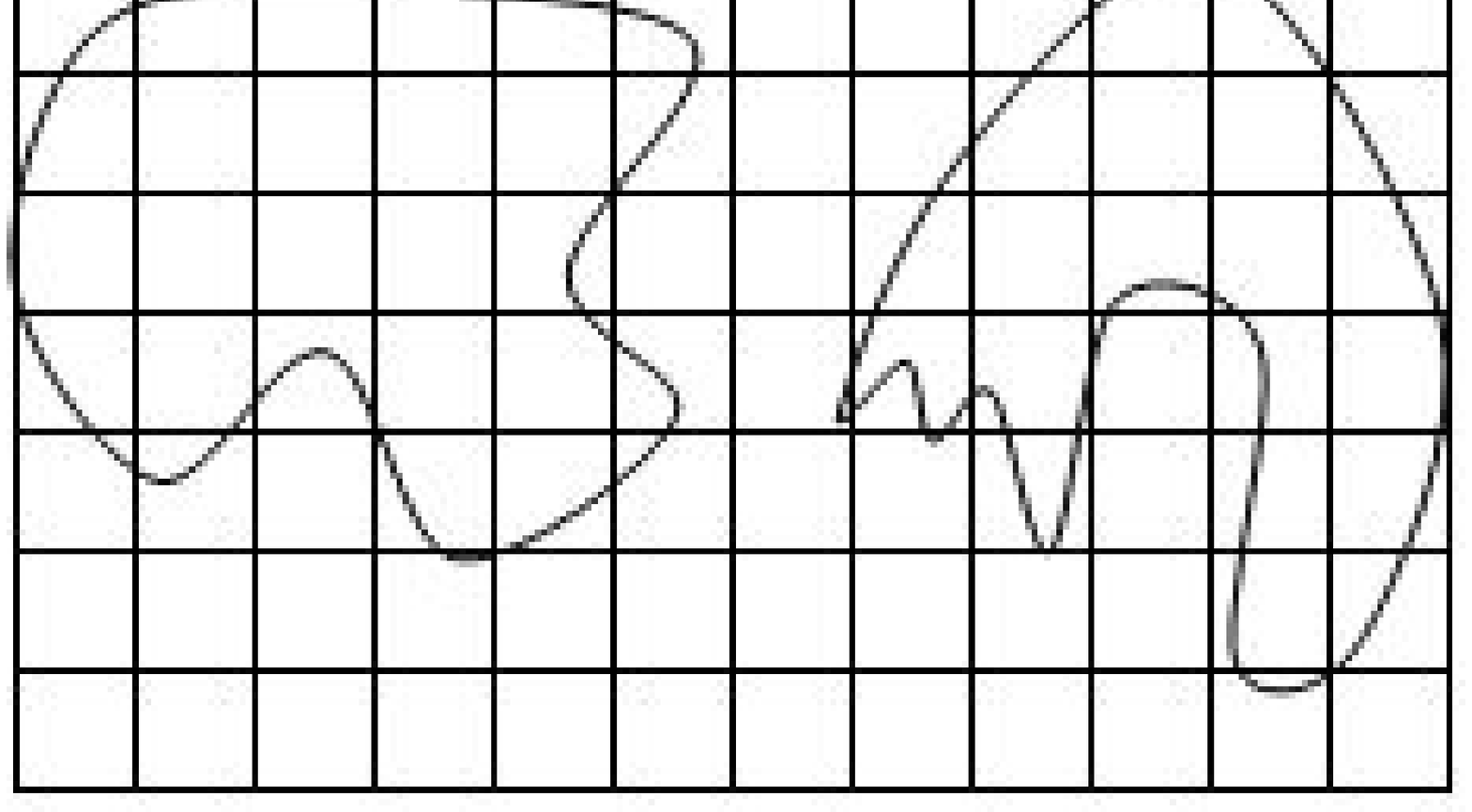}\includegraphics[width=0.3\linewidth]{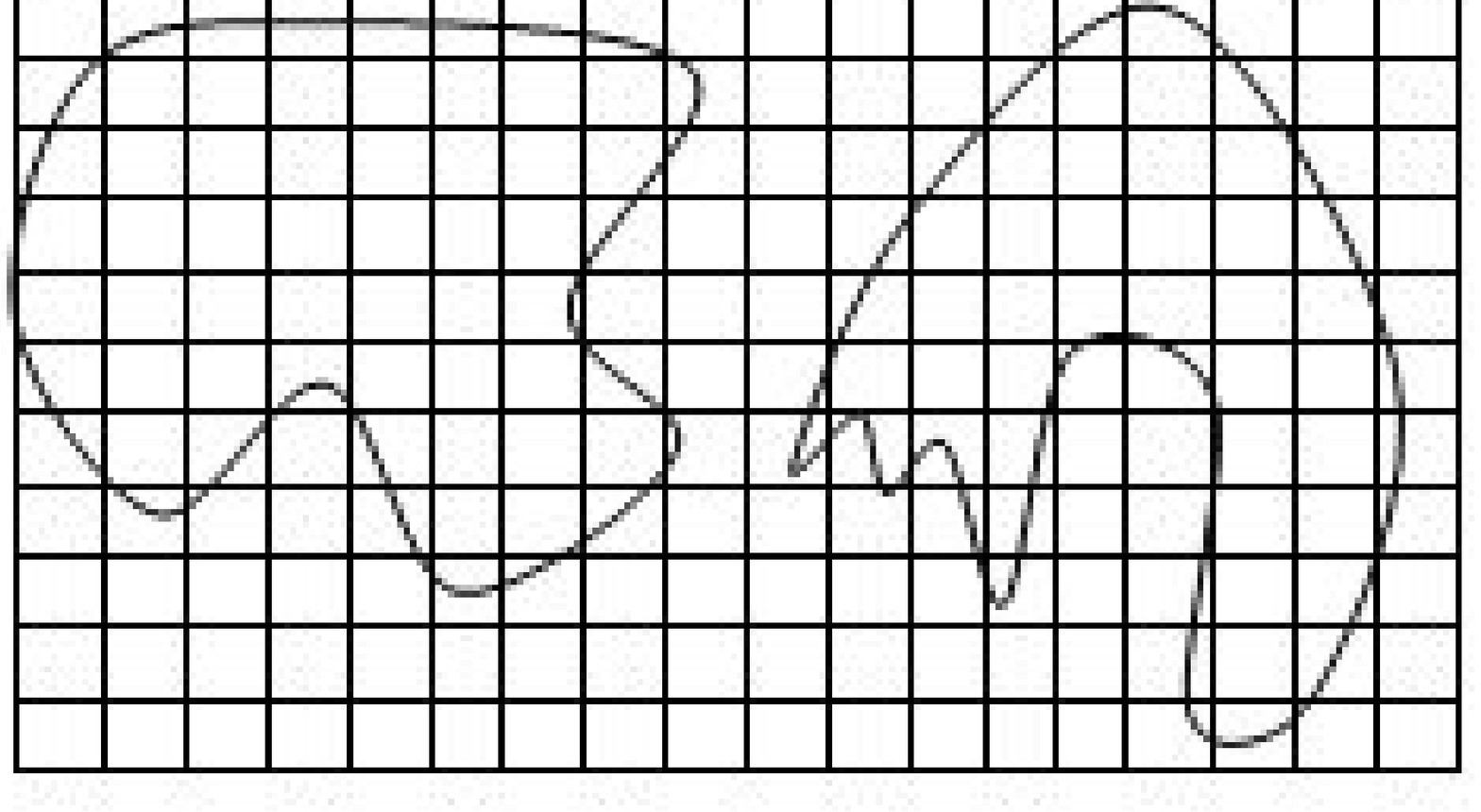}
		\vskip-4cm
		\caption{Area of planeidom}
		\label{fig:p2}
	\end{figure}
}

You can already see that the smaller ,,scale'' enables more accurate work. The length, area, and volume of simple geometric shapes already in antiquity
they were known and countable. The concept of area was extended by Peano and Jordan
for a larger system of subsets of the plane in the XIX. at the end of the century.

First, we will talk about the subsets of the space $\R^2$. By a square of the space $\R^2$ we always mean a rectangle whose sides are parallel to the axes $X_1X_2$, that is, $\R^2$ has a

\begin{equation} \label{first}
T=\{(x_1,x_2)\in\R^2\vert a \leq x_1\leq b,c\leq x_2 \leq d\}
\end{equation}
shaped subset. The area of this

\[
t=(b-a)(d-c).
\]
The interval in \eqref{first} can be divided as follows:

\[
\begin{split}
a
&=
x_1^{(0)}<\ldots<x_1^{(n)}=b\\
c
&=
x_2^{(0)}<\ldots<x_2^{(m)}=d
\end{split}
\]
unequally, then $x_1=x_1^{(i)}$ $(i=0,1,\ldots,n)$, and $x_2=x_2^{(k)}$ $(k=0,1,\ldots,m)$ is called a square grid, the subset calld to elementary quadrilaterals or lattice angles.

\[
T_{ik}=\big\{(x_1,x_2)\in\R^2 \vert x_1^{(i-1)}\leq x_1\leq x_1^{(i)}, x_2^{(k-1)}\leq x_2 \leq x_2^{(k)}\big\}.
\]

In the following, let $A\subset\R^2$ be a bounded set. Then there exists a quadrilateral $T$ such that $A\subset T$. In this case, we can say that $T$ covers the set $A$. In the same sense, we can also say that the set $A$ is also covered by the quadrilateral grid obtained by dividing the sides of the quadrilateral $T$, since it is obvious that $A\subset\cup_{i,k}=T$, where the union is all $T$ refers to its elementary quadrilateral.

Now we examine the set of quadrilaterals covering $A$ and the sum of the areas of those elementary quadrilaterals that contain at least one point of $A$. The lower bound of these area sums is called the outer area of the set $A$. And the internal area of the set $A$ is the upper limit of the sum of the areas of the elementary quadrilaterals inside the set $A$, provided that such quadrilaterals exist at all. If there are none, then the internal area of $A$ is considered to be zero.

\subsection{Jordan measure}

After this, we can interpret the concept of measurability and area of the set $A$. We say that the set $A\subset\R^2$ has Jordan-measurable area (or Jordan-measurable for short) if the outer area and inner area of $A$ are equal. The common measure of the outer area and the inner area is called the area of the set $A$.

According to this, the exterior measure of a planar bounded set according to Jordan should be the finite covering it
the exact lower bound of the area of many polygonal shapes, according to Jordan
and the internal measure of the area of the finitely many polygonal figures lying within it
exact upper limit. If these are equal, then the set is Jordan-measurable, and the common value is called the Jordan measure of the set. So we are used to this measure
in primary and secondary schools to be called an area. The above connections are illustrated in the following Figure \ref{fig:jordan}.

\begin{figure} [h!]
	\centering
	\includegraphics[width=0.5\linewidth]{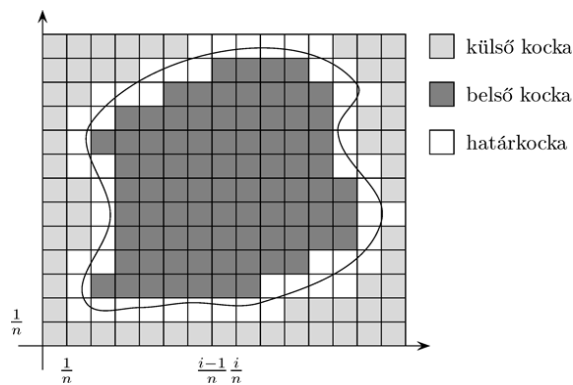}
	\caption{Covering outer and inner cubes}
	\label{fig:jordan}
\end{figure}

The square below was divided into small squares with side lengths of $1/n$, and in each case we marked which part is the outer part or the inner part of the given planar block, by definition the boundary cube belongs to both parts, so it cannot be marked exactly. The larger $n$ we take, the more accurate value we get. The Jordan measure and the Riemann integral are very closely related, since they are one
a non-negative real function can be integrated exactly according to Riemann if the plane under the curve of the function is Jordan-measurable. Then the Riemann integral and a
the Jordan measurement of the planeidom the same.

The previously sketched procedure for marking the area of subsets in $\R^2$ can also be used to interpret the volume of bounded subsets of the space $\R^3$, and an analogous procedure can also be used in the case of length measurement.

In three-dimensional space, of course

\[
T=\big\{\textbf{x}\in\R^3\vert \textbf{x}=(x_1,x_2,x_3), a_i\leq x_i \leq b_i, i=1,2,3 \big \}
\]
we work with three-dimensional bricks of the shape, the volume of which is a
\[
A=(b_1-a_1)(b_2-a_2)(b_3-a_3)
\]
is defined by the formula, instead of square grids, we use so-called brick grids. The outer and inner volume of a bounded set $A\subset\R^3$ can be interpreted with the help of such a brick grid in the same way as we saw in $\R^2$, and the Jordan measurability and volume measure are also the same as the concepts seen previously.

\begin{fa}
	\label{circle}
	Determine the area of a circle with radius $r(>0)$!
\end{fa}

\ans{The idea here is similar to the methods seen above, let's approximate the area from two sides using regular polygons.
	
	\begin{figure}[h!]
		\centering
		\includegraphics[width=0.25\linewidth]{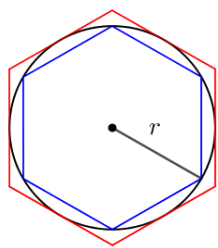}
		\caption{Approximation of the circumference of a circle for $n=6$}
		\label{fig:korkozelites}
	\end{figure}
	
	\[
	\begin{split}
		A_{\text{inscribed polygon}}
		&\leq
		A_{\text{circle}}\leq A_{\text{polygon inscribed around it}}\\
		\dfrac{r^2\sin\left(\frac{2\pi}{n}\right)}{2}\cdot n
		&\leq
	    A_{\text{circle}}\leq\dfrac{2r^2\text{tg}\left(\frac{\pi}{n}\right)}{2}\cdot n,
	\end{split}
	\] 
but 
\[\lim_{n\to\infty}\dfrac{r^2\sin\left(\frac{2\pi}{n}\right)}{2}\cdot n=\lim_{n\to\infty}\dfrac{r^2}{2}\cdot 2\pi \cdot \dfrac{\sin\left(\frac{2\pi}{n}\right)}{\frac{2\pi}{n}}=r^2\pi,\] 
since
\[\lim_{x\to0}\dfrac{\sin x}{x}=1.\] 
Furthermore
\[\lim_{n\to \infty} r^2\cdot \text{tg}\left(\frac{\pi}{n}\right)\cdot n=\lim_{n\to \infty}r^2\pi \dfrac{1}{\cos\left(\frac{\pi}{n}\right)}\cdot \dfrac{\sin\left(\frac{\pi}{n}\right)}{\frac{\pi}{n}}=r^2\pi.\]
Then, based on the previously mentioned police principle, we get that 
\[r^2\pi\leq A_{\text{circle}}\leq r^2\pi,\]
so $A_{\text{circle}}=r^2\pi$.
}

\begin{fa}
	Determine the volume of the straight circular cylinder with base radius $r(>0)$ and height $m$.
\end{fa}

\ans{To determine the volume of a straight circular cylinder, we will use the fact that the volume of the column is given as follows:
\begin{equation} \label{eq:hasab}
	V_{\text{prism}}=A_{\text{surface area}}\cdot m_{\text{prism}}.
\end{equation}}

\begin{figure}[h!]
	\centering
	\includegraphics[width=0.4\linewidth]{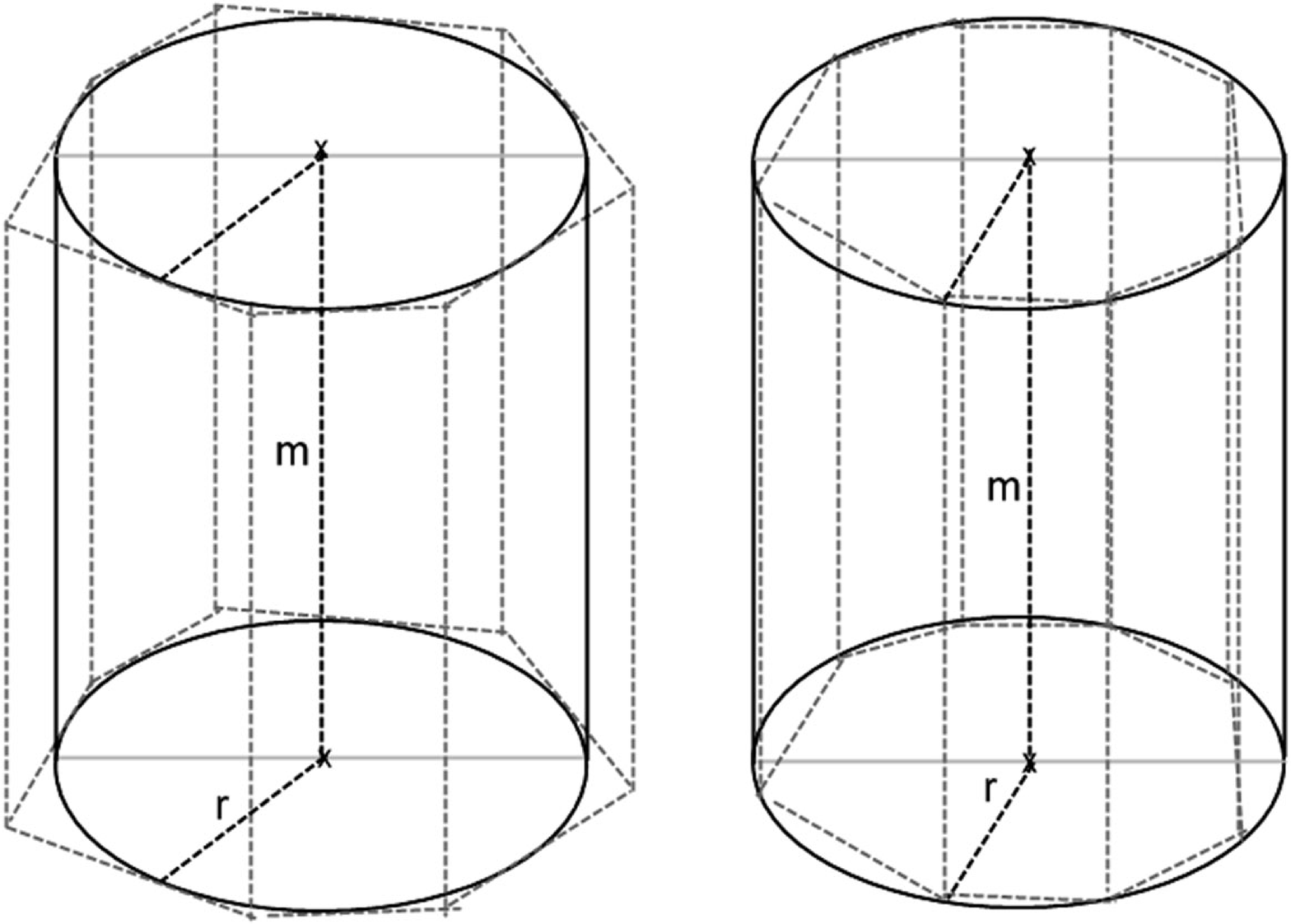}
	\caption{Approximation of the cylinder with a column}
	\label{fig:hengera}
\end{figure}

The volume of the straight cylinder is determined by the two-sided approximation method using columns inscribed around it, as shown in Figure \ref{fig:hengera}. 

The volume of a straight cylinder is determined by the two-sided approximation method using columns inscribed around it and inside it. In the base of the cylinder, i.e., in the circle with radius $r$, and around the circle, we write a regular polygon, the number of sides of which is denoted by $n$. In the case of inscribed columns, the vertices of the polygon are located on the outline, and in those written around it, the sides of the polygons (the base edges of the base and cover plate) touch the base of the cylinder or cover circle. The base and top plates of the columns and the cylinder fall in one plane. The volume of the inscribed columns is always smaller, and the volume of the surrounding columns is always greater than the volume of the cylinder, so we can write the following inequalities:

\begin{equation}
	A_{\text{into}}\cdot m=V_{\text{into it}}\leq V_{\text{cylinder}}\leq V_{\text{around}}=A_{\text{around}}\cdot m,
\end{equation}
where $m$ is the height of the cylinder and columns, $A_{\text{into}}$ and $A_{\text{around}}$ are the areas of the bases of the columns written in and around them, respectively. By increasing the number of sides of the inscribed and circumscribed columns, the area of the inscribed polygon and thus the volume of the inscribed column increases, while the area of the circumscribed polygon and the volume of the column decreases. By increasing the number of sides of the polygons inscribed in and around it, the difference between the areas of the two polygons can be made as small as possible, and this gives the area of the circle, i.e. $r^2\pi$ (see \textbf{Exercise \ref{circle}.}). Thus, the volume of the cylinder: $V_{\text{cylinder}}=r^2\pi m$.

\subsection{Synthetic geometry - Proof of elementary theorems: The Pythagorean theorem}

\begin{figure}[h!]
		\centering
		\includegraphics[width=0.30\linewidth]{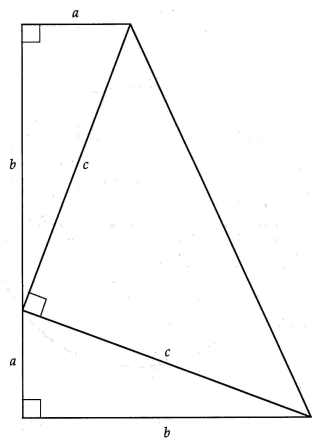}
		\caption{Proving the Pythagorean Theorem}
		\label{fig:graf}
\end{figure}

Among the many proofs of the theorem, one of the simplest is based on writing the area of a right-angled trapezoid in two ways. Indeed, writing the area of the trapezoid alone or as the sum of its parts: $$(a+b)^2=\frac{2ab+c^2}{2}.$$ Rearranging this, $$a^2+b^2= c^2$$ results\footnote{This proof is said to be attributed to the 20th President of the United States, James A. Garfield.} (Figure \ref{fig:graf}). 
	
\subsection{Synthetic geometry - Proving compatibility theorems: Reversal of Ceva's theorem}

Let us first consider the following concept and a related theorem.

\begin{def*}
	
In a triangle, it connects the vertices to a point on the opposite side sections are called Ceva\footnote{Giovanni Ceva, Italian mathematician (1648-1736), 1678.} sections if they intersect at one point (Figure \ref{ceva1}).
	
	\begin{figure}[h!]
		\centering
		\includegraphics[width=0.30\linewidth]{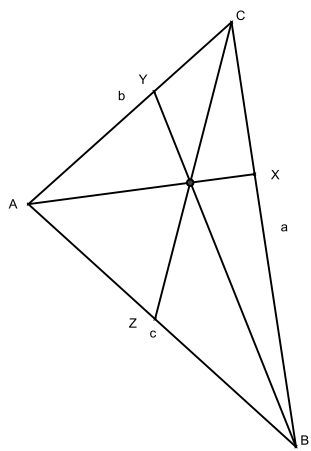}
		\caption{Ceva's sections}
		\label{ceva1}
	\end{figure}
	
\end{def*} 

\begin{theorem}
	In the triangle $ABC$, for the Ceva segments $AX, BY, CZ$ we got $$  \frac{BX}{XC}\frac{CY}{YA}\frac{AZ}{ZB}=1.$$
\end{theorem}

\begin{fa}
 Prove the inversion of the Ceva theorem.
\end{fa}

\begin{theorem}[The inversion of the Ceva theorem]
	If in the triangle $ABC$ for the segments $AX, BY, CZ$ where $X, Y, Z$ are respectively the interior points of the sides $BC, CA, AB$
\end{theorem}

\begin{equation}\label{cv1}
\frac{BX}{XC}\frac{CY}{YA}\frac{AZ}{ZB}=1
\end{equation}
exists, then the three sections fit into one point.

,,Zoom in'' to the intersection point(s) in two ways, as follows. Denote the intersection of $AX$ and $BY$ by $P$. Denote the end point of the third Ceva section passing through $P$ by $Z'$. Then, due to Ceva's theorem, \begin{equation}\label{cv2}
\frac{BX}{XC}\frac{CY}{YA}\frac{AZ'}{Z'B}=1
 \end{equation} and due to (\ref{cv1}), $$\frac{BX}{XC}\frac{CY}{YA}\frac{AZ}{ZB}=1,$$ so from (\ref{cv1}) and from (\ref{cv2}) we get: $$\frac{AZ'}{Z'B}=\frac{AZ}{ZB}\Longrightarrow Z=Z'.$$
	
\subsection{Synthetic geometry - Solving fitting problems}

Below we show two examples of a sharper appearance of our method.

\begin{fa}
	We connected the vertices $AE$ and $DF$ of the squares $ABCD$ and $BEFG$ shown in the figure. We prove that the connecting sections intersect on the $BG$ side (Figure \ref{fig:negyzet1})! \begin{center}\begin{figure}[h!]
			\centering
			\includegraphics[width=0.3\linewidth]{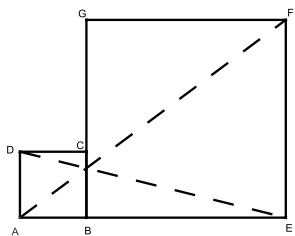}
			\caption{The squares and connecting sections}
			\label{fig:negyzet1}
		\end{figure}
		
	\end{center}
\end{fa}

\ans{Denote the intersection of $AF$ and $GB$ segments by $H,$ and let the intersection of $DE$ and $GB$ segments be $I.$ Also let $HB=x$ and $IB=y$ (Figure \ref{fig:negyzet2}). We must realize that $H=I,$ and thanks to our new notations, this is equivalent to $x=y.$ Let the side length of the square $ABCD$ be $a$, and that of the square BEFG $b$. 
	
\begin{center}
		\begin{figure}[h!]
			\centering
			\includegraphics[width=0.30\linewidth]{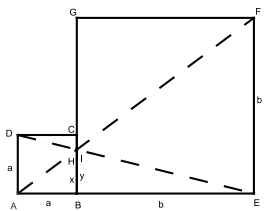}
			\caption{Squares and intersections}
			\label{fig:negyzet2}
		\end{figure}
\end{center}

Note that $AHB\triangle\sim AFE\triangle$ because their sides are pairwise parallel. Following from the similarity, the ratio of the lengths of the corresponding sides is the same, i.e. $\frac{x}{a}=\frac{b}{a+b}$. Similarly, $EIB\triangle\sim EDA\triangle$, so $\frac{y}{b}=\frac{a}{a+b}$. Then, comparing the equations, we get: $x=y$, i.e. $H=I$, i.e. the two sections really intersect at the $BG$ section.
}

\begin{fa}
Let the inscribed circle of the triangle $ABC$ be the point of contact on the side $c$ of the triangle $E.$ Prove that the inscribed circles of the triangle $AEC$ and the triangle $CEB$ touch each other (Figure \ref{fig:haromszog1}) . \begin{center}
		\begin{figure}[h!]
		\centering
		\includegraphics[width=0.35\linewidth]{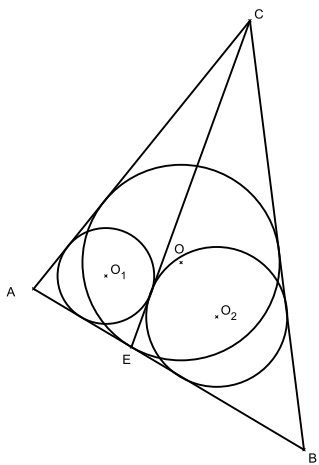}
		\caption{The triangle and the inscribed circles}
		\label{fig:haromszog1}
	\end{figure}
	\end{center}
	\end{fa}

\ans{Let the inscribed circle of the triangle $ABC$ be $k$, its center $O$, the inscribed circle of the triangle $ACE$ be $k_1$, its center $O_1$, also denote the inscribed circle of the triangle $CEB$ by $k_2$, the center of the circle $ O_2.$ Let the points of contact of circles $k_1$ and $k_2$ with $CE$ be $D$ and $F$, respectively (Figure \ref{fig:haromszog2}).
	\begin{center}
		\begin{figure}[h!]
			\centering
			\includegraphics[width=0.35\linewidth]{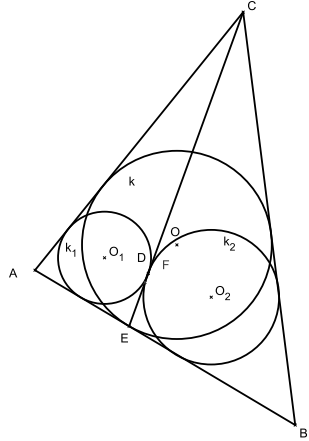}
			\caption{The triangle and the points of contact}
			\label{fig:haromszog2}
		\end{figure}
	\end{center}
	Then, based on the summation of the tangent points of the inscribed circle for the lengths of the segments cut from the sides of the triangle from the triangle $ACE$:
	$$DE=s_1-b=\frac{CE+b+s-a}{2}-b=\frac{CE+s-b-a}{2},$$ and from the triangle $CEB$: $$FE=s_2-a=\frac{CE+a+s-b}{2}-a=\frac{CE+s-b-a}{2}$$
	 can be written where $s_1$ is $ACE$ and $s_2$ is the semicircumference of triangle $CEB$. Now $DE=FE$ or $E=F$ results from comparing the equations.
}

\section{Probability}

Using our previous results and experience for infinite series, we can successfully apply our method to the subject of probability calculation.

\subsection{A probability model, a geometric series}

There may be tasks where the probability of certain events can only be described using infinite series. Using another model, however, we can get results without knowing infinite series, and we can also arrive at a method for summing infinite series from two different approaches to the given problem.

\begin{fa}
Two players, Anna and Balázs, play with a regular dice. Anna starts the game and the dice are thrown alternately one at a time, and the first to roll a six wins. What is the probability that Anna wins?
\end{fa}

\ans{Let $A$ denote the event that Anna wins. Let's find the value of $\p(A)$. Let us introduce the notation $\p(A_n)$, which means that Anna wins in step $n$. We know that the events $\{A_n\}_{n=1}^\infty$ are disjoint, so using the $\sigma$-additive property of probability, $\p(A)=\p \left(\cup_{ n=1}^\infty A_n\right)=\sum_{n=1}^\infty \p(A_n)$. It is easy to see that $n$ can only be an odd number, since Anna can only win on the odd number of throws. Let $n=2k+1$ $(k \in \mathbb{N})$, then $\p(A_1)=\frac{1}{6}$, $\p(A_3)=\frac{1 }{6}(\frac{5}{6})^2$, since in order for Anna to win on the 3rd throw, the first throw she throws and the second throw thrown by Balázs must not be a six, while Anna the 3rd roll thrown by should be exactly a six. Similar considerations result in $\p(A_n)=\frac{1}{6}(\frac{5}{6})^{n-1}$. Thus, we can write the sought probability in the form of the following sum: $$\sum_{n=1}^\infty \p(A_n)=\sum_{k=0}^\infty \p(A_{2k+1})= \sum_{k=0}^\infty \frac{1}{6}\left(\frac{5}{6}\right)^{2k}.$$
	Now let's see another approach. Consider the following figure (Figure \ref{fig:valsz1}). 
	\begin{center}
		\begin{figure}[h!]
			\centering
			\includegraphics[width=0.3\linewidth]{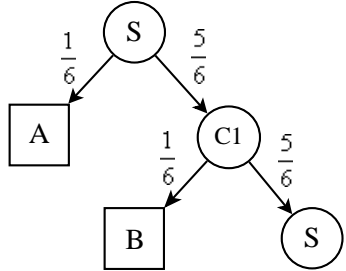}
			\caption{States and their probabilities}
			\label{fig:valsz1}
		\end{figure}
	\end{center}

Here, $S$ denotes the start state, hence the player A, who throws a six with probability $\frac{1}{6}$, thus we get to the winning state A, or does not throw a six with probability $\frac{5}{6}$ , so we get to state $C_1$, from where B throws. If you roll a six, the probability of this is still $\frac{1}{6}$, and if you do not roll a six, then A again follows, so we have returned to state $S$. Now we are looking for the probability of winning A from the start state, denote this by $\p_S(A)$, for the other states as well. Then we can arrive at the following system of equations.

\[
\begin{split}
\p_S(A)
&=
\dfrac{1}{6}+\dfrac{5}{6}\p_{C_{1}}(A)\\
\p_{C_{1}}(A)
&=
\dfrac{5}{6}\p_S(A)
\end{split}
\]
from which $P_S(A)=\frac{6}{11}$ results. Thus, the probability of the same event was calculated in two ways:
	
\[
\sum_{k=0}^\infty \frac{1}{6}\left(\frac{5}{6}\right)^{2k}=\frac{6}{11}.
\]
}

\section{Some more difficult tasks}

To illustrate the greatness of the method, here are some tasks that are much more difficult than before.

\subsection{Fibonacci numbers and non-Fibonacci numbers}

\begin{fa}
	Prove that if $0 < m < n$, then $X=f_{2m+1}+f_{2m+3}+ \ldots +f_{2n+1}$ and $Y=f_{2m }+f_{2m+2}+ \ldots + f_{2n}$ numbers are not members of the Fibonacci sequence, where $f_n$ denotes the $n$th Fibonacci number.
\end{fa}

\ans{We will prove that $X$ and $Y$ fall between two adjacent Fibonacci numbers, so they cannot be Fibonacci numbers. Consider the expressions $X=(X+Y)-Y$ and $Y=X-(X-Y)$. Then, $$X=f_{2m}+f_{2m+1}+f_{2m+2}+f_{2m+3}+ \ldots f_{2n}+f_{2n+1}-f_{2m}-f_{2m+2}- \ldots - f_{2n}=$$ $$=f_{2m+2}+f_{2m+4}+ \ldots + f_{2n+2}-f_{2m}-f_{2m+2}- \ldots - f_{2n}=f_{2n+2}-f_{2m}\hskip 5pt \hbox{and}$$ $$Y=f_{2m+1}+f_{2m+3}+ \ldots + f_{2n+1}-(f_{2m+1}-f_{2m}+f_{2m+3}-f_{2m+2}+ \ldots + f_{2n+1}-f_{2n})=$$$$=f_{2m+1}+f_{2m+3}+ \ldots + f_{2n+1}-(f_{2m-1}+f_{2m+1} \ldots f_{2n-1})=f_{2n+1}-f_{2m-1}.$$ Since $m < n$, so $f_{2m} < f_{2n}$, and $f_{2m-1} < f_{2n-1}$, so $$f_{2n+2} > X = f_{2n+2} - f_{2m} >f_{2n+2} - f_{2n} = f_{2n+1}$$ and $$f_{2n+1} > Y = f_{2n+1}-f_{2m-1} > f_{2n+1}- f_{2n-1}=f_{2n}$$ i.e. both X and Y fall between two adjacent Fibonacci numbers, and thus cannot be a Fibonacci number.}

\subsection{A probability model, a non-geometric series}

\begin{fa} 
Let us prove that $$\frac{1}{2}\cdot\left(1+\frac{(-1)^{n+1}}{3^n}\right)=\sum_{L=1}^\infty\frac{\binom{2L}{n-1}}{2^{2L+1}}$$
\end{fa}

\ans{Consider the following problem and solve it in two ways.
	
\begin{fa} Anna and Balázs play by tossing a regular coin. The winner is the one who throws the $n$th head. What is the probability that Anna wins if she starts the game?
\end{fa}}

\end{document}